%% file: manuscript-amspreprint.tex
\pdfoutput=1

\IfFileExists{./prepreamble-amspreprint.sty}{\RequirePackage[packages,theorems,changes]{./prepreamble-amspreprint}}{}

\makeatletter
\IfFileExists{./scoop-latex/scoop-packages.sty}{\providecommand*{\input@path}{}\edef\input@path{{./scoop-latex/}\input@path}}{}
\@namedef{ver@minted.sty}{}
\@namedef{opt@minted.sty}{}
\makeatother
\PassOptionsToPackage{style=authoryear-comp,backend=biber}{biblatex}
\PassOptionsToPackage{commentmarkup = footnote}{changes}    

\documentclass[english]{amsart}

\RequirePackage{biblatex}
\RequirePackage{ltxcmds}
\IfFileExists{./preamble-amspreprint.sty}{\RequirePackage[packages,theorems,changes]{./preamble-amspreprint}}{}

\usepackage{scoop-local}

\IfFileExists{./postpreamble-amspreprint.sty}{\RequirePackage[packages,theorems,changes]{./postpreamble-amspreprint}}{}

\makeatletter
\@ifpackageloaded{changes}{
\definechangesauthor[name = {Frederik Köhne}, color = {red!80!black}]{FK}
\definechangesauthor[name = {Leonie Kreis}, color = {blue!80!black}]{LK}
\definechangesauthor[name = {Anton Schiela}, color = {orange!80!black}]{AS}
\definechangesauthor[name = {Roland Herzog}, color = {green!70!black}]{RH}
}{}
\makeatother        

\makeatletter
\@ifpackageloaded{hyperref}{%
	\hypersetup{
		pdftitle = {Adaptive Step Sizes for Stochastic Gradient Descent in Hilbert Spaces},
		pdfauthor = {Frederik Köhne, Leonie Kreis, Anton Schiela, Roland Herzog},
		pdfkeywords = {stochastic gradient descent, adaptive learning rates, preconditioning, strongly convex problems},
		pdfcreator = {Created using the Scoop Template Engine version 1.3.0.}
	}
}{
	\pdfinfo{
		/Title (Adaptive Step Sizes for Stochastic Gradient Descent in Hilbert Spaces)
		/Author (Frederik Köhne, Leonie Kreis, Anton Schiela, Roland Herzog)
		/Subject ()
		/Keywords (stochastic gradient descent, adaptive learning rates, preconditioning, strongly convex problems)
		/Creator (Created using the Scoop Template Engine version 1.3.0.)
	}
}
\makeatother

\title[Adaptive Step Sizes for SGD in Hilbert Spaces]{Adaptive Step Sizes for Stochastic Gradient Descent in Hilbert Spaces}

\author[F. Köhne]{Frederik Köhne\orcidlink{0009-0008-6185-9675}}
\address[F. Köhne]{Department of Mathematics, University of Bayreuth, 95440 Bayreuth, Germany}
\email{\detokenize{frederik.koehne@uni-bayreuth.de}}
\urladdr{https://num.math.uni-bayreuth.de/en/team/frederik-koehne/}

\author[L. Kreis]{Leonie Kreis\orcidlink{0000-0003-4234-4867}}
\address[L. Kreis]{Interdisciplinary Center for Scientific Computing, Heidelberg University, 69120 Heidelberg, Germany}
\email{\detokenize{leonie.kreis@iwr.uni-heidelberg.de}}
\urladdr{https://scoop.iwr.uni-heidelberg.de}

\author[A. Schiela]{Anton Schiela\orcidlink{0000-0002-6959-2951}}
\address[A. Schiela]{Department of Mathematics, University of Bayreuth, 95440 Bayreuth, Germany}
\email{\detokenize{anton.schiela@uni-bayreuth.de}}
\urladdr{https://num.math.uni-bayreuth.de/en/team/anton-schiela/}

\author[R. Herzog]{Roland Herzog\orcidlink{0000-0003-2164-6575}}
\address[R. Herzog]{Interdisciplinary Center for Scientific Computing, Heidelberg University, 69120 Heidelberg, Germany}
\address[R. Herzog]{Institute for Mathematics, Heidelberg University, 69120 Heidelberg, Germany}
\email{\detokenize{roland.herzog@iwr.uni-heidelberg.de}}
\urladdr{https://scoop.iwr.uni-heidelberg.de}

\thanks{This work was supported by DFG grants HE~6077/13--1 and SCHI~1379/8--1 within the Priority Program SPP~2298 (Mathematical Foundations of Deep Learning), which is gratefully acknowledged.}

\date{\today}

\dedicatory{}

\begin{document}

\begin{abstract}
\input{abstract.tex}
\end{abstract}


\makeatletter
\ltx@ifpackageloaded{hyperref}{%
\subjclass[2010]{\href{https://mathscinet.ams.org/msc/msc2020.html?t=65K05}{65K05}, \href{https://mathscinet.ams.org/msc/msc2020.html?t=68T05}{68T05}, \href{https://mathscinet.ams.org/msc/msc2020.html?t=68Q25}{68Q25}, \href{https://mathscinet.ams.org/msc/msc2020.html?t=90C15}{90C15}}
}{%
}
\makeatother

\maketitle

\input{content.tex}

\appendix
\input{appendix.tex}

\printbibliography

\end{document}

%% file: abstract.tex
This paper proposes a novel approach to adaptive step sizes in stochastic gradient descent (SGD) by utilizing quantities that we have identified as numerically traceable --- the Lipschitz constant for gradients and a concept of the local variance in search directions.
Our findings yield a nearly hyperparameter-free algorithm for stochastic optimization, which has provable convergence properties and exhibits truly problem adaptive behavior on classical image classification tasks.
Our framework is set in a general Hilbert space and thus enables the potential inclusion of a preconditioner through the choice of the inner product.

%% file: content.tex
\section{Introduction}
\label{section:introduction}

Stochastic Gradient Descent (SGD) is a simple, yet effective algorithm commonly used to solve stochastic optimization problems.
These problems only allow access to a noisy, usually unbiased estimate of the target function, called the \emph{sampled function}, and its derivative in each iteration.
Formally, one is interested in minimizing
\begin{equation*}
	F(w)
	=
	\int_\Omega f_\xi(w) \d P(\xi)
	,
\end{equation*}
where $P$ is a probability measure on $\Omega$
and for every $\xi \in \Omega$, $f_\xi \colon \R^n \to \R$ is a suitable function. In each iteration $k$, $\xi_k \sim P$ is sampled and $f_{\xi_k}$ is used as the sampled function.
Such problems prominently arise in modern machine learning applications, where selecting $f_{\xi_k}$ corresponds to the selection of one or few addends in a finite sum optimization problem.
SGD, first introduced in \cite{RobbinsMonro:1951:1}, has since become the workhorse for this kind of problems and has led to the development of several variants of the algorithm.

\subsection{Known Adaptive Step Size Strategies}
\label{subsection:known-adaptive-step-size-strategies}

It is well known that the performance of SGD, as well as convergence guarantees, crucially depend on the step sizes (\emph{learning rates}) employed.
Therefore, different approaches to make the step sizes of SGD adaptive have been developed.
We briefly summarize them in what follows.

\subsubsection{Polyak-Type Strategies}

Polyak-type strategies aim to adapt the well known Polyak step sizes for classical gradient schemes, first presented in~\cite{Polyak:1987:1}, to the stochastic setting.
A common assumption is that the minimum, or at least a lower bound to the minimum of the sampled function, is known.
In \cite{LoizouVaswaniLaradjiLacosteJulien:2021:1}, the authors derive convergence properties of SGD with Polyak-type step sizes for the interpolating setting (no noise at the minimizer) and convergence to a neighborhood of the minimizer for the non-interpolating case (noise present at the minimizer); see \cref{subsection:noise-at-minimizer} for a discussion of the two settings.
\cite{JiangStich:2023:1} extended the work of \cite{LoizouVaswaniLaradjiLacosteJulien:2021:1} to obtain convergence, even in the non-interpolating setting.
For strongly convex target functions, both works obtain linear convergence in the interpolating setting.
The latter work also shows sublinear convergence of order $O(\frac{1}{\sqrt k})$ in the non-interpolating setting.

\subsubsection{Line Search Strategies}

Line search strategies aim to apply the concept of line search from classical optimization to stochastic optimization.
These strategies typically involve the repeated evaluation of the sampled function at various candidates for the next iterate until a desired decrease is observed.
A direct adaptation of the well-established Armijo line search is documented in \cite{VaswaniMishkinLaradjiSchmidtGidelLacosteJulien:2019:1}.
Convergence theory for line search methods must often consider the noise introduced by the sampled function.
A theory that controls this noise can be found in \cite{PaquetteScheinberg:2020:1}.
Both works achieve linear convergence in the strongly convex, interpolating regime.
A significant limitation of line search methods is the repeated evaluation of the sampled function at each iteration, which can become computationally expensive.

\subsubsection{Diagonal Scaling Methods}

Another class of commonly used adaptive methods can be classified as diagonal scaling methods, which gather information from past iterations to develop a step size strategy where each dimension of the input space has its unique step size.
It is also possible to interpret these as methods that employ a \emph{diagonal} preconditioning matrix to the derivative, in order to obtain the search direction, where the preconditioning matrix needs not to be constant over time.
Prominent examples of these methods include RMSProp\footnote{Proposed in unpublished work \cite{HintonSrivastavaSwersky:2012:1} by Geoffrey Hinton et al.; see also \cite{Ruder:2016:1}.}, Adagrad (\cite{DuchiHazanSinger:2011:1}), Adadelta (\cite{Zeiler:2012:1}), and Adam (\cite{KingmaBa:2015:1}), as well as its numerous variants.
For these classes of algorithms it remains, however, unclear how the choices of scalings are related the convergence of the algorithms.
In \cite{VaswaniLaradjiKunstnerMengSchmidtLacosteJulien:2020:1} the authors propose to use line search methods to set up the step size for Adagrad.

\subsubsection{Trust Region Methods}

Another line of research focuses on trust region methods.
Here, adaptivity stems from selecting the trust region radius based on previous iterations.
Examples of such work can be found in \cite{BlanchetCartisMenickellyScheinberg:2016:1} and \cite{CurtisShi:2020:1}.

\subsection{Variance in the Search Direction}
\label{subsection:variance-in-the-search-direction}

The primary theoretical concern when examining stochastic optimization methods such as SGD is the variance in the search direction, and specifically, the variance in the gradient of the sampled function.
Consequently, methods should either control the variance (\ie, be able to reduce it to an appropriate threshold), or manage it in real time without being able to control it directly.
Gradient aggregation techniques such as SVRG (\cite{JohnsonZhang:2013:1}) and SAGA (\cite{DefazioBachLacosteJulien:2014:1}) incorporate supplementary gradients or gradients from preceding iterations to decrease noise.
This approach typically necessitates extensive computational and/or memory resources.
Algorithms aiming to manage noise typically do so under special assumptions regarding the variance of the search direction.
These assumptions involve a globally bounded second moment of the search direction, as outlined in~\cite{NemirovskiJuditskyLanShapiro:2009:1, HazanKale:2014:1}.
More recently, a variance bound was proposed by \cite{BottouCurtisNocedal:2018:1}, which allows for non-zero variance at stationary points and growth of the variance proportional to the squared norm of the true gradient.
These bounds can be used to select a step size.
In~\cite{NguyenNguyenVanDijkRichtarikScheinbergTakac:2018:2}, it was shown by the authors that the assumption in \cite{BottouCurtisNocedal:2018:1} is dispensable and can be deduced from specific smoothness and convexity assumptions on the problem.

\subsection{Noise at the Minimizer}
\label{subsection:noise-at-minimizer}

The behavior of SGD is significantly influenced by the noise at the minimizer $w^\star$, which can be measured by $\E{\xi}{\norm{\nabla f_\xi(w^\star)}^2}$.
Two regimes are distinguished in the literature that lead to qualitatively completely different behavior of SGD.
The first regime, the \emph{interpolating} setting, corresponds to the absence of noise at the minimizer.
In this case, $w^\star$ is not only a critical point for~$F$, but also for all~$f_\xi$.
In machine learning, this case corresponds to the heavily overparameterized case, where the model is capable of interpolating the training data.

The \emph{non-interpolating} setting corresponds to the presence of noise at the minimizer.
In this case, there are sampled functions that are not stationary at~$w^\star$.
Thus, even if the algorithm arrives at $w^\star$, it will not recognize this and will even move away.

Generally speaking, the interpolating setting allows much stronger convergence results (\cite{JiangStich:2023:1,MaBassilyBelkin:2018:1,GarrigosGower:2023:1}).
For instance, the interpolating setting allows for convergence with a constant, positive step size, while this does not hold for the non-interpolating case.
In \cite{JiangStich:2023:1}, the authors argue that adaptive step size schemes should be \emph{robust} in the sense that they automatically adapt to the setting (interpolating \versus non-interpolating).

\subsection{Our Contribution}
\label{subsection:contribution}

In this paper we propose adaptive step size schemes, based on key quantities governing the convergence behavior of SGD.
These quantities describe, on the one hand, the nonlinearity of the problem, and, on the other hand, its stochasticity.
Our step size schemes use computable estimates for these quantities to control the progress of the iteration.
We analyze these schemes in the framework of $\mu$-strongly convex and $L$-smooth functions.

We present a step size scheme that are not affected by the strong convexity parameter~$\mu$.
We argue that this is crucial, as $\mu$ is usually neither available, nor can it be estimated reasonably.
By contrast, the smoothness constant~$L$ can in fact be estimated, as outlined in \cref{section:estimation-techniques}.
It is well known in the literature that step sizes proportional to the inverse of the smoothness constant~$L$ are known to make progress initially, but fail to converge to the optimum in the non-interpolating setting (\cite{GarrigosGower:2023:1}).
Our adaptive scheme, presented in \cref{section:convergence-to-optimality}, resolves this issue by incorporating an additional estimator for the local variance.
We show how the variance can be estimated by reusing information that was already captured during the estimation of~$L$.
With this step size strategy, we obtain a robust adaptive step size in the sense of \cref{subsection:noise-at-minimizer}.
We show linear convergence in the interpolating case and convergence of order $1/k$ in the non-interpolating case.

\subsection{Outline}
\label{subsection:outline}

We first review existing and commonly used variance models and their implications on adaptive step size selection.
Motivated by the fact that these models are a consequence of smoothness and convexity (\cite{NguyenNguyenVanDijkRichtarikScheinbergTakac:2018:2}), we investigate the influence of convexity.
In \cref{proposition:unbounded-variation-bound-P,proposition:unbounded-variation-bound-Pstar} we show that the constants in the variance bounds can become arbitrarily large when certain properties of the family~$f_\xi$ can't be controlled, or when the strong convexity parameter decreases to zero.
We conclude that these variance models are not suited for adaptive step size selection because they cause the step size to depend on the strong convexity parameter.

We proceed by examining modified versions of the variance bound, which lead to descent results for SGD with step sizes independent of the strong convexity parameter (\cref{lemma:adapted-variance-bound,theorem:descent-SGD}).

In \cref{section:convergence-to-optimality} we identify a step size strategy that only invokes quantities that can be estimated numerically.
We give a convergence-in-expectation proof for SGD with this step size strategy.
The result (\cref{theorem:SGD-global-convergence}) shows linear convergence in the interpolating regime and convergence of order $O(\frac{1}{k})$ in the non-interpolating regime.

In \cref{section:estimation-techniques}, we describe methods capable of estimating the quantities invoked by the step size strategy.
For classical neural network applications, the estimation process costs one additional forward pass at each mini-batch.
In \cref{section:practical-aspects} we comment on the details of implementation for our method, which pitfalls we expect and how safeguards against them could be established.
In \cref{section:numerical-results} we finally show how our method performs, on synthetic problems, as well as on classical image classification tasks.

\section{Problem Setting}
\label{section:problem-setting}

Our results do not require the Euclidean structure of~$\R^n$, nor are they confined to finite dimensions.
We therefore consider a real Hilbert space~$X$ as our setting.
The inner product is denoted by $\skp{\cdot}{\cdot}$.
The dual space of~$X$ is denoted by~$X^*$, and the dual pairing between $X$ and $X^*$ is written as $\dual{\cdot}{\cdot}$.
The Riesz isomorphism is $H \colon X \to X^*$, which maps $x \mapsto \skp{\cdot}{x} \in X^*$.
The derivative~$f'$ of a differentiable function $f \colon X \to \R$ is a mapping $f' \colon X \to X^*$.
The gradient is the Riesz representer of the derivative, \ie, $\nabla f(w) = \inv{H} f'(w)$.
In other words, $\nabla f(w) \in X$ is the unique element that satisfies $\skp{\nabla f(w)}{v} = f'(w) \, v$ for all $v \in X$.

(Stochastic) gradient methods on $X$ compute (an estimate of) $f'(w)$ and a search direction $\delta w = -\inv{H} f'(w)$.
This highlights that selecting the inner product of~$X$ --- thus choosing the Riesz isomorphism --- corresponds to selecting a preconditioner for iterative, gradient-based methods.
In the case $X = \R^n$, equipped with the standard inner product, we obtain the classical (stochastic) gradient method with $H = \id_n$.

Suppose that $(\Omega, \cA, P)$ is a probability space\footnote{The $\sigma$-algebra $\cA$ will not play a significant role in our work.}
Suppose that $h \colon \Omega \to \R, \xi \mapsto h(\xi)$ is a measurable function.
The quantity
\begin{equation*}
	\E{\xi}{h}
	=
	\int_\Omega h(\xi)\d P(\xi)
\end{equation*}
is the expected value or mean of $h$.
Let $\xi_0, \ldots, \xi_k$ be independent realizations of the random variable $\xi \sim P$.
For the measurable function
\begin{equation*}
	h
	\colon
	\Omega^{k-s+1}
	\to
	\R
	,
	\;
	(\xi_s, \dots, \xi_k)
	\mapsto
	h(\xi_s, \dots, \xi_k)
\end{equation*}
with $0 \le s \le k$, we define the short-hand notation
\begin{equation*}
	\E{k:k}{h}
	\coloneqq
	\E{\xi_k}h
	\quad
	\text{and recursively}
	\quad
	\E{s:k}{h}
	\coloneqq
	\E[auto]{\xi_s}{ \E{s+1:k}{h}}
	.
\end{equation*}
In this work we consider stochastic optimization problems of the following form:
\begin{definition}
	\label{definition:SOP}
	Suppose that $f_\xi \colon X \to \R$ is differentiable for each $\xi \in \Omega$, and that $\xi \mapsto f_\xi(w)$ is integrable for every $w \in X$.
	Define
	\begin{equation*}
		F(w)
		\coloneqq
		\int_\Omega f_\xi(w) \d P(\xi)
		=
		\E{\xi}{f_\xi(w)}
		.
	\end{equation*}
	We refer to the problem of finding $w^\star \in X$ with $F(w^\star) = \inf_{w \in X} F(w)$ as the \emph{stochastic optimization problem} (SOP) with data $(f_\xi, \Omega, P)$, or short the SOP $(f_\xi, \Omega, P)$.
\end{definition}

\begin{definition}
	\label{definition:L-smoothness-mu-strong-convexity}
	Suppose that $f \colon X \to \R$ is a differentiable function with derivative $f' \colon X \to X^*$.
	\begin{itemize}
		\item
			$f$ is said to be $L$-smooth for some $L \ge 0$ if
			\begin{equation*}
				\norm{f'(w_1) - f'(w_2)}_{X^*}
				\le
				L \, \norm{w_1 - w_2}_X
			\end{equation*}
			holds for all $w_1, w_2 \in X$.

		\item
			$f$ is said to be $\mu$-strongly convex for some $\mu \ge 0$ if
			\begin{equation*}
				\dual{f'(w_1) - f'(w_2)}{w_1 - w_2}
				\ge
				\mu \, \norm{w_1 - w_2}_X^2
			\end{equation*}
			holds for all $w_1, w_2 \in X$.
	\end{itemize}
\end{definition}

\begin{definition}
	\label{definition:mu-L-feasibility}
	Suppose that $0 < \mu \le L < \infty$.
	\begin{enumerate}
		\item
			A differentiable mapping $f \colon X \to \R$ is said to be $(\mu, L)$-feasible if $f$ is $\mu$-strongly convex and $L$-smooth.

		\item
			An SOP $(f_\xi, \Omega, P)$ according to \cref{definition:SOP} is said to be $(\mu, L)$-feasible if its mean $F$ is $(\mu, L)$-feasible.

		\item
			An SOP $(f_\xi, \Omega, P)$ is said to be strongly $(\mu, L)$-feasible if $f_\xi$ is $(\mu, L)$-feasible for almost every $\xi \in \Omega $.

		\item \label{item:pointwise-feasible}
			Suppose $\xi \mapsto \mu_\xi, \xi \mapsto L_\xi \colon \Omega \to [0, \infty)$ are measurable functions.
			An SOP $(f_\xi, \Omega, P)$ is said to be pointwise $(\mu_\xi, L_\xi)$-feasible if $f_\xi$ is $(\mu_\xi, L_\xi)$-feasible for almost every $\xi \in \Omega$ and $\Lmax \coloneqq \esssup_{\xi \in \Omega} L_\xi < \infty$.
	\end{enumerate}
\end{definition}

\begin{remark}
	In \cref{item:pointwise-feasible} of the preceding definition we explicitly include the case $\mu_\xi = 0$ for \emph{some} $\xi \in \Omega$.
	For such values of~$\xi$, the function $f_\xi$ is convex but not strongly convex.
\end{remark}

By the Lax-Milgram Theorem, each $(\mu,L)$-feasible SOP has a unique minimizer~$w^\star$.

We denote the \emph{variance} of the noisy derivative $f_\xi'(w)$ by%
\footnote{%
	In the classical setting $X = \R^n$ with dual space $X^* = \R^{1 \times n}$ and Riesz isomorphism $H v = v^\transp$, we have
	\begin{equation*}
		\Var{\xi}{\nabla f_\xi(w)}
		=
		\Var{\xi}{f_\xi'(w)}
		=
		\E[auto]{\xi}{\norm{\nabla f_\xi(w) - \nabla F(w)}_2^2}
		.
	\end{equation*}
}
\begin{equation*}
	\Var{\xi}{f_\xi'(w)}
	\coloneqq
	\E[auto]{\xi}{\norm{f_\xi'(w) - F'(w)}_{X^*}^2}
	.
\end{equation*}
Note that this is not precisely the variance of a vector valued quantity as usually defined in statistics, but rather the trace of its covariance matrix.
However, it enjoys the same separation property,
\begin{equation}
	\Var{\xi}{f_\xi'(w)}
	=
	\E[auto]{\xi}{\norm{f_\xi'(w)}_{X^*}^2} - \norm{F'(w)}_{X^*}^2
	.
	\label{equation:var-shift-property}
\end{equation}

\subsection{SGD Descent Analysis}
\label{subsection:SGD-descent-analysis}

A straightforward and therefore popular algorithm for solving a stochastic optimization problem is Stochastic Gradient Descent (SGD), first proposed by~\cite{RobbinsMonro:1951:1}.
In each iteration, an independent realization of a random variable $\xi$ with distribution $P$ is computed and then $\nabla f_\xi(w)$ is used as a search direction, which is an unbiased estimator of $\nabla F(w)$.
SGD then performs the simple update step $w^+ = w - \alpha \nabla f_\xi(w)$, where $\alpha > 0$ is a step size or \emph{learning rate}, see \cref{algorithm:sgd-basic}.
\begin{algorithm}[htp]
	\caption{SGD}
	\label{algorithm:sgd-basic}
	\begin{algorithmic}[1]
		\Require Step size $\alpha > 0$, initial iterate $w_0$.
		\For{$k \ge 0$}
		\State Sample $\xi_k \sim P$.
		\State $w_{k+1} \gets w_k - \alpha \nabla f_{\xi_k}(w_k)$.
		\EndFor
	\end{algorithmic}
\end{algorithm}
An overview of descent properties of SGD algorithms can be found in \cite{BottouCurtisNocedal:2018:1}.
Their results are based on a popular assumption of the type
\begin{equation}
	\label{eq:var_assumption}
	\Var{\xi}{\nabla f_\xi(w)}
	\le
	V_0 + V_1 \, \norm{\nabla F(w)}_X^2
	\quad
	\text{for all }
	w \in X
	,
\end{equation}
which describes the variance at the minimizer with some ground noise $V_0$ and allows for growth of the variance further away from the minimizer.
Using a constant step size $\alpha$, the authors in \cite{BottouCurtisNocedal:2018:1} show linear convergence in expectation of the suboptimality gap $F(w_k) - F({w^\star})$ to a stagnation level $\frac{\alpha V_0 L}{2 \, \mu}$:
\begin{equation}
	\E[auto]{}{F(w_{k+1}) - F(w^\star) - \frac{\alpha V_0 L}{2 \, \mu}}
	\le
	(1 - \alpha \, \mu) \, \paren[auto](){F(w_k) - F(w^\star) - \frac{\alpha V_0 L}{2 \, \mu}}
	,
	\label{eq:convergence-BottouCurtisNocedal:2018:1}
\end{equation}
where the expectation is taken over the randomness in the search directions.
This result holds for any \emph{sufficiently small} step size, namely
\begin{equation}
	\alpha
	\le
	\frac{1}{L \, (1 + V_1)}
	.
	\label{eq:bound-on-alpha-L-V1}
\end{equation}

\subsection{Problems Arising}
\label{subsection:problems-arising}

In order to establish an adaptive scheme for the step size, based on the model described above, one needs access to the quantities that affect the step size in \eqref{eq:bound-on-alpha-L-V1}, namely $L$ and $V_1$.
As these are usually unknown in practice, it is necessary to find estimators, which use numerically traceable quantities.
The following lemma shows that $V_1$ can be bounded in terms of $L_\xi$ and $\mu$.

\begin{restatable}{lemma}{lemmafirstboundonvariance}
	\label{lemma:first-bound-on-variance}
	Let $(f_\xi, \Omega, P)$ be a $(\mu, L)$-feasible SOP such that $f_\xi$ is $L_\xi$-smooth for some measurable function $\xi \mapsto L_\xi$.
	Then the variance assumption \eqref{eq:var_assumption} holds with
	\begin{equation*}
		V_0
		=
		2 \, \E[auto]{\xi}{\norm{\nabla f_\xi(w^\star)}_X^2}
		\quad
		\text{and}
		\quad
		V_1
		=
		2 \, \frac{ \E{\xi}{L_\xi^2}}{\mu^2} - 1
		.
	\end{equation*}
\end{restatable}
\begin{proof}
	The proof can be found in \cref{section:proof-for-variance-behavior}.
\end{proof}

Under the stronger assumption $L_\xi \le \Lmax$ almost everywhere, we even obtain:
\begin{restatable}{lemma}{lemmastrongerboundonvarianceuniformlyboundedLcase}
	\label{lemma:stronger-bound-on-variance-uniformly-bounded-L-case}
	Let $(f_\xi, \Omega, P)$ be a pointwise $(\mu_\xi, L_\xi)$-feasible SOP such that $F$ is $\mu$-strongly convex for some $\mu > 0$.
	Then the variance assumption \eqref{eq:var_assumption} holds with
	\begin{equation*}
		V_0
		=
		2 \, \E[auto]{\xi}{\norm{\nabla f_\xi(w^\star)}_X^2}
		\quad
		\text{and}
		\quad
		V_1
		=
		2 \, \frac{\Lmax}{\mu} - 1
		.
	\end{equation*}
\end{restatable}
\begin{proof}
	The proof can be found in \cref{section:proof-for-variance-behavior}.
\end{proof}

Thus, if $L_\xi$ is uniformly bounded, we obtain an improved bound for the variance compared to \cref{lemma:first-bound-on-variance}.
A similar analysis can be found in \cite{NguyenNguyenVanDijkRichtarikScheinbergTakac:2018:2}.

\subsection{Asymptotic Behavior of the Variance}
\label{subsection:asymptotic-behavior-of-the-variance}

The bounds on $V_1$ established in the previous section depend on~$\mu$ and exhibit the property $V_1 \to \infty$ when $\mu \to 0$.
Thus, when using a step size $\alpha = \frac{1}{L \, (1 + V_1)}$, as suggested by \eqref{eq:bound-on-alpha-L-V1}, this implies~$\alpha \to 0$ when $\mu \to 0$, which is clearly not desirable.
In the following we exhibit that this behavior is a necessary consequence of the structure of the bound \eqref{eq:var_assumption}.
Our results show that it is impossible to select~$V_0$ and~$V_1$ in a way that \eqref{eq:var_assumption} holds and the constants are not affected by~$\mu$ or~$\Var{\xi}{L_\xi}$.

\begin{definition}
	Given an SOP $(f_\xi, \Omega, P)$ following \cref{definition:SOP} and $V_0 \ge 0$, let
	\begin{equation*}
		V_1(V_0)
		\coloneqq
		\sup \setDef[auto]{\frac{\Var{\xi}{\nabla f_\xi(w)} - V_0}{\norm{\nabla F(w)}_X^2}}{\Var{\xi}{\nabla f_\xi(w)} > V_0, \; w \ne w^\star}
	\end{equation*}
	denote the smallest possible constant~$V_1$ such that the variance assumption \eqref{eq:var_assumption} is met.
\end{definition}

\begin{restatable}{proposition}{propositionunboundedvariationboundP}
	\label{proposition:unbounded-variation-bound-P}
	Suppose that $\cP(\mu, L)$ is the set of all $(\mu, L)$-feasible stochastic optimization problems $(f_\xi,\Omega,P)$.
	Then for any $\mu \in (0,1)$ we have
	\begin{equation*}
		\sup_{(f_\xi,\Omega,P) \in \cP(\mu, 1)}
		\inf_{V_0 \in \R} V_1(V_0)
		=
		\infty
		.
	\end{equation*}
\end{restatable}
\begin{proof}
	The proof can be found in \cref{section:proof-for-variance-behavior}.
\end{proof}
Thus, the constant~$V_1$ in bounds of the type of \eqref{eq:var_assumption} can become arbitrarily large for certain distributions.
In the proof of \cref{proposition:unbounded-variation-bound-P}, which is given in \cref{section:proof-for-variance-behavior}, we used a heavy-tailed distribution to let $\Var{\xi}{L_\xi}$ grow arbitrarily, which leads to the variance of the gradient growing arbitrarily, while $\nabla F(w)$ remains bounded.

Such behavior can not occur if we consider strongly $(\mu,L)$-feasible problems.
This can be seen easily by observing that $\Var{\xi}{L_\xi} = \E[big]{\xi}{\abs{L_\xi - \E{\xi}{L_\xi}}^2} < L^2$, since $\E{\xi}{L_\xi} \in [\mu, L]$ and therefore $\abs{L - \E{\xi}{L_\xi}} < L$.
Recalling \cref{lemma:stronger-bound-on-variance-uniformly-bounded-L-case}, we see that for strongly-$(\mu,L)$-feasible problems, we obtain the stronger bound $V_1 \le 2 \, \frac{L}{\mu} - 1$.
However, this still becomes arbitrarily large with $\mu \to 0$.
This also is not a flaw in the results, but rather a necessary consequence, as the following result demonstrates.

\begin{restatable}{proposition}{propositionunboundedvariationboundPstar}
	\label{proposition:unbounded-variation-bound-Pstar}
	Suppose that $\cP_\star(\mu, L)$ is the set of all strongly-$(\mu, L)$-feasible stochastic optimization problems $(f_\xi,\Omega,P)$.
	Then for any $\mu \le \frac{1}{2}$ we have
	\begin{equation*}
		\sup_{(f_\xi,\Omega,P) \in \cP_\star(\mu, 1)}
		\inf_{V_0 \in \R} V_1(V_0)
		\ge
		\frac{1}{64 \mu}
		.
	\end{equation*}
	In particular, the bound in \cref{lemma:stronger-bound-on-variance-uniformly-bounded-L-case} is asymptotically sharp.
\end{restatable}
This result is also proved in \cref{section:proof-for-variance-behavior}.
In the proof, we provide a family of strongly-$(\mu, L)$-feasible SOPs, whose selection is geometrically motivated; see \cref{remark:geometric-motivation-for-variance-result}.
\makeatother
We use the same class of SOPs in the experiment presented in \cref{figure:steps-size-need-not-depend-on-mu}.

\subsection{\texorpdfstring{Variance Bounds Independent of $\mu$}{Variance Bounds Independent of Convexity}}
\label{subsection:variance-bounds-and-initial-convergence}

\begin{figure}[ht]
	\centering{\includegraphics[width = 1\textwidth]{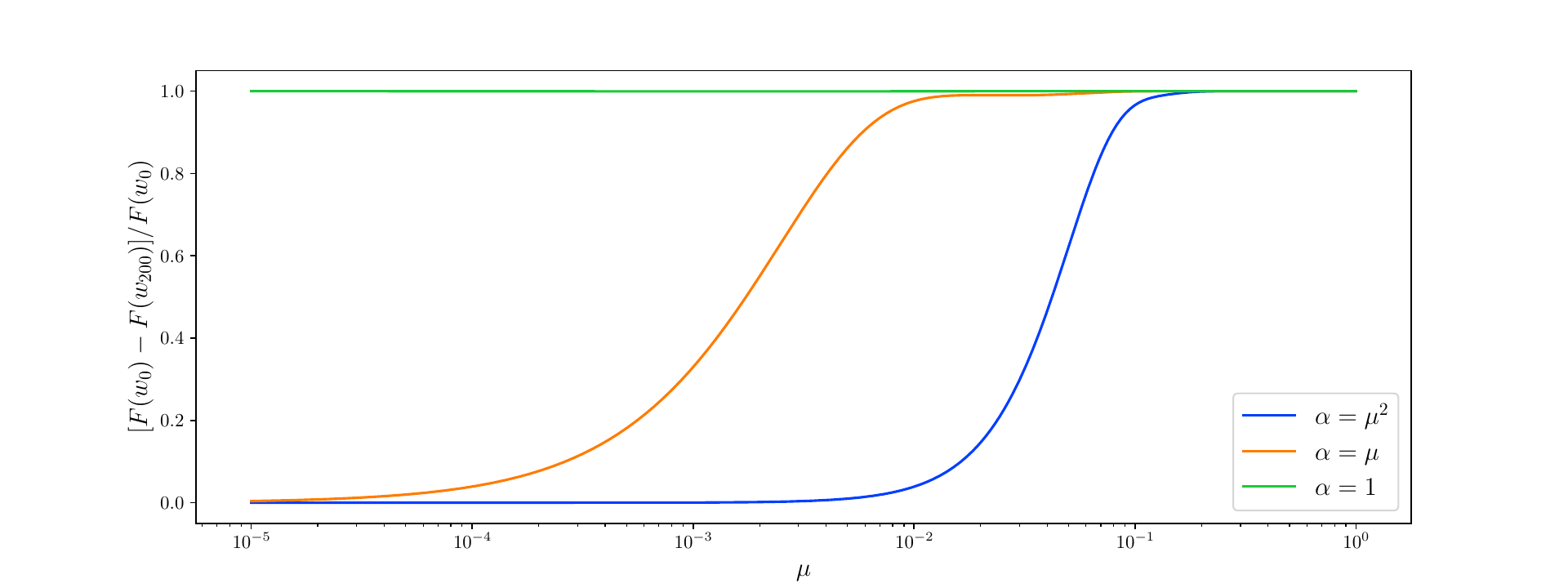}}
	\caption{A step size $\sim \mu$ is too conservative.
		The figure shows a comparison of different step sizes, in dependency of the convexity parameter~$\mu$ for the example in the proof of \cref{proposition:unbounded-variation-bound-Pstar}.
		SGD's relative progress is plotted, with higher values indicating better performance.
		According to the theory presented in \cref{section:problem-setting}, a step size of $\frac{1}{L \, (1 + V_1)}$ should be employed.
		As shown in the proof of \cref{proposition:unbounded-variation-bound-Pstar}, $V_1$ grows at a rate of $\frac{1}{\mu}$ in this example.
		Therefore, keeping $L = 1$ fixed would result in a step size of $\sim \mu$.
	However, this approach appears to be too conservative.}
	\label{figure:steps-size-need-not-depend-on-mu}
\end{figure}

\Cref{proposition:unbounded-variation-bound-P,proposition:unbounded-variation-bound-Pstar} show that the constants in a variance bound of the form \eqref{eq:var_assumption} can become arbitrarily large even for seemingly harmless problems.
The behavior implied by \cref{proposition:unbounded-variation-bound-Pstar} in particular is undesirable, since it leads to an unnecessary reduction in the step size in case the latter is chosen according to \eqref{eq:bound-on-alpha-L-V1}, as \cref{figure:steps-size-need-not-depend-on-mu} illustrates.

A step size that decreases with $\mu$ is also not well suited for adaptive methods because $\mu$ is difficult to estimate.
For problems with quadratic objective, $\mu$ is the smallest eigenvalue of the Hessian matrix, which is basically as hard to estimate as estimating all eigenvalues.
Additionally, the matrix itself is not available.
Instead, only noisy results of matrix-vector products are computed in this case.
This easily leads to the conclusion that neither can knowledge of~$\mu$ be expected, nor is the estimation of~$\mu$ possible in practical scenarios.

Consequently, a model for the variance as in~\eqref{eq:var_assumption} is not well suited for determining step sizes.
In the following, we present a different model for the variance, which allows for bounds provably independent of $\mu$.
Results in this direction are already present in the literature.
The main idea is to replace $\norm{F'(w)}_{X^*}^2$ in \eqref{eq:var_assumption} by another quantity, which only scales \emph{linearly} in the convexity constant $\mu$.
If $w$ is chosen appropriately, $\norm{F'(w)}_{X^*}^2$ is proportional to $\mu^2$, which explains the dependence of $V_1$ on $\mu$ in the results above.\\
For the remainder of this section denote
\begin{equation*}
	V_0
	=
	\E[auto]{\xi}{\norm{f_\xi'(w^\star)}_{X^*}^2}
\end{equation*}
and assume $V_0 < \infty$.
\begin{lemma}[Adapted Variance Bound]
	\label{lemma:adapted-variance-bound}
	Suppose that $(f_\xi, \Omega, P)$ is a pointwise $(\mu_\xi, L_\xi)$-feasible SOP.
	For $w \in X$, denote $D_w = F(w) - F(w^\star)$.
	It holds:
	\begin{align*}
		\E[auto]{\xi}{\norm{f_\xi'(w)}_{X^*}^2}
		&
		\le
		4 \, \Lmax \, D_w
		+
		2 \, V_0
		\\
		&
		\le
		4 \, \Lmax \, \dual{F'(w)}{w - w^\star}
		+
		2 \, V_0
	\end{align*}
	and thus
	\begin{align*}
		\Var[big]{\xi}{f_\xi'(w)}
		&
		\le
		4 \, \Lmax \, D_w
		+
		2 \, V_0
		-
		\norm{F'(w)}_{X^*}^2
		\\
		&
		\le
		4 \, \Lmax \, \dual{F'(w)}{w - w^\star}
		+
		2 \, V_0
		-
		\norm{F'(w)}_{X^*}^2
		.
	\end{align*}
\end{lemma}
\begin{proof}
	The bound
	\begin{equation*}
		\E[auto]{\xi}{\norm{f_\xi'(w)}_{X^*}^2}
		\le
		4 \, \Lmax \, D_w
		+
		2 \, V_0
	\end{equation*}
	can be found in \cite[Lemma~4.20]{GarrigosGower:2023:1}, see also \cite{JohnsonZhang:2013:1}.
	In these references, the authors consider the finite sum setting and the case $X = \R^d$, equipped with the standard euclidean inner product.
	As their arguments can be applied directly to our setting, we refrain from providing a more detailed proof here.
	Due to convexity we have $D_w \le \dual{F'(w)}{w - w^\star}$, which implies the second bound.
	The bound on the variance follows from $\Var[big]{\xi}{f_\xi'(w)} = \E[big]{\xi}{\norm{f_\xi'(w)}_{X^*}^2} - \norm{F'(w)}_{X^*}^2$.
\end{proof}

With the bounds from \cref{lemma:adapted-variance-bound} in place, the following convergence to stagnation result can be shown:

\begin{theorem}
	\label{theorem:descent-SGD}
	Suppose that $(f_\xi, \Omega, P)$ is a pointwise $(\mu_\xi, L_\xi)$-feasible SOP such that $F$ is $\mu$-strongly convex.
	Denote by $w_k$ the sequence generated by SGD with a constant step size $0 < \alpha \le \frac{1}{2\Lmax}$.
	Then
	\begin{equation*}
		\E[auto]{0:k-1}{\norm{w_k - w^\star}_X^2}
		\le
		(1 - \mu \, \alpha)^k \, \norm{w_0 - w^\star}_X^2
		+
		2 \, \frac{\alpha V_0}{\mu}
		.
	\end{equation*}
\end{theorem}
\begin{proof}
	This result can be found in \cite[Theorem~5.8]{GarrigosGower:2023:1}.
	As in the case of \cref{lemma:adapted-variance-bound}, the proofs can easily be adapted to our setting.
\end{proof}

\section{Adaptive Step Sizes and Convergence to Optimality}
\label{section:convergence-to-optimality}

The results presented in the previous sections only ensure convergence to a stagnation level, since the step size suggested there ($\sim \frac{1}{\Lmax}$) does not tend to zero closer to the minimizer.
In this section, we devise a method to determine step sizes that do not only depend on the nonlinearity of the problem ($\Lmax$), but also take into account the local noise.
The suggested step sizes are proven to lead to convergence of order at least $O(1/k)$ in the non-interpolating setting and to linear convergence in the interpolating setting.
Therefore, the method matches the best known convergence rates in the respective settings.

From $w_{k+1} = w_k - \alpha_k \nabla f_\xi(w_k) = w_k - \alpha_k \inv{H} f_{\xi_k}'(w_k)$ we infer
\begin{align}
	\E{\xi_k}{F(w_{k+1})}
	&
	\le
	\E[auto]{%
	\xi_k}{F(w_k)
		-
		\alpha_k \, \dual{F'(w_k)}{f_{\xi_k}'(w_k)}
		+
		\frac{L \, \alpha_k^2}{2} \, \norm{f_{\xi_k}'(w_k)}_{X^*}^2
	}
	\notag
	\\
	&
	=
	F(w_k)
	-
	\alpha_k \norm{F'(w_k)}_{X^*}^2
	+
	\frac{L \, \alpha_k^2}{2} \, \E[big]{\xi_k}{\norm{f_{\xi_k}'(w_k)}_{X^*}^2}
	.
	\label{eq:descent-alpha}
\end{align}
Minimizing the right-hand side \wrt $\alpha$ suggests using the step size
\begin{equation}
	\label{eq:perfect-step-size}
	\alpha_k
	=
	\frac{\norm{F'(w_k)}_{X^*}^2}{L \, \E[big]{\xi_k}{\norm{f_{\xi_k}'(w_k)}_{X^*}^2}}
	=
	\frac{ \E[big]{\xi_k}{\norm{f_{\xi_k}'(w_k)}_{X^*}^2} - \Var{\xi_k}{f_{\xi_k}'(w_k)}}{L \, \E[big]{\xi_k}{\norm{f_{\xi_k}'(w_k)}_{X^*}^2}}
	.
\end{equation}
Here, we used the shift property \eqref{equation:var-shift-property}.
We argue that all quantities which are used in \eqref{eq:perfect-step-size} are traceable during the execution of SGD.
We refer to \cref{subsection:nonlinearity-estimation} for the estimation of the nonlinearity ($L$), to \cref{subsection:variance-estimation} for the variance, and to \cref{subsection:accessibility-of-expected-value-gradient-norm} for $\E{\xi_k}{\norm{f_{\xi_k}'(w_k)}_{X^*}^2}$.

The step size discussed in \cref{section:problem-setting} only incorporates nonlinearity and a bound on the asymptotic growth of the variance, which describes the behavior of the variance for large $\norm{w_k - w^\star}$.
It is thus not surprising that it can only lead to convergence to a stagnation level, which is determined by the noise at the minimizer.
In contrast to that, the step size suggested in~\eqref{eq:perfect-step-size} incorporates the true local noise.
This allows for convergence to the minimizer, even in the presence of noise at the minimizer, \ie, in the non-interpolating case.
\begin{remark}
	The step sizes in~\eqref{eq:perfect-step-size} can also be written as
	\begin{equation*}
		\alpha_k
		=
		\frac{1}{L}\paren[auto](){1 - \frac{\Var{\xi_k}{f_{\xi_k}'(w_k)}}{\E[auto]{\xi_k}{\norm{f_{\xi_k}'(w_k)}_{X^*}^2}}}
		.
	\end{equation*}
	This highlights that two factors determine good step sizes.
	On the one hand, this is the nonlinearity of the problem, described by $L$.
	On the other hand, we have the factor $\paren[auto](){1 - \frac{\Var{\xi_k}{f_{\xi_k}'(w_k)}}{\E[b g]{\xi_k}{\norm{f_{\xi_k}'(w_k)}_{X^*}^2}}} \in [0,1]$, describing the local stochasticity of the problem.
	Recall that $\E[auto]{\xi_k}{\norm{f_{\xi_k}'(w_k)}_{X^*}^2} = \norm{F'(w_k)}_{X^*}^2 + \Var{\xi_k}{f_{\xi_k}'(w_k)}$.
	Thus, if the variance is relatively small compared to $\norm{F'(w_k)}_{X^*}^2$, we have step sizes close to $\frac{1}{L}$, as we would expect in the deterministic setting.
	Conversely, if the variance becomes relatively large compared to $\norm{F'(w_k)}_{X^*}^2$, as it would be the case close to a minimizer in the non-interpolating setting, the step sizes also become small.
	Such behavior of the step sizes is not needed in deterministic optimization, but is crucial in the non-interpolating setting in stochastic optimization.
\end{remark}

\subsection{Convergence Analysis}
\label{subsection:convergence-analysis}

Inserting the step sizes from \eqref{eq:perfect-step-size} into \eqref{eq:descent-alpha}, we obtain the following chain of inequalities, denoting $D_k \coloneqq F(w_k) - F(w^\star)$:
\begin{align}
	\E[auto]{\xi_k}{D_{k+1}}
	&
	\le
	D_k -\frac{\alpha_k}{2} \norm{F'(w_k)}_{X^*}^2
	\notag
	\\
	&
	\le
	(1 - \mu \, \alpha_k) \, D_k
	\quad
	\text{(due to strong convexity)}
	\label{eq:linear-convergence-mu-alpha}
	\\
	&
	\le
	D_k \exp\paren[auto](){-\mu \, \alpha_k}
	\quad
	(\text{due to $1 - x \le \exp(x)$})
	\notag
	\\
	&
	\le
	D_k \exp\paren[bigg](){-2 \, \frac{\mu^2}{L}\frac{D_k}{\E[big]{\xi_k}{\norm{f_{\xi_k}'(w_k)}_{X^*}^2}}}
	.
	\notag
\end{align}
In the last inequality, we used \eqref{eq:perfect-step-size} and $\norm{F'(w_k)}_{X^*}^2 \ge 2 \, \mu D_k$; this argument was used in~\eqref{eq:linear-convergence-mu-alpha} as well.
To state our convergence theorem, we need the following two lemmas.

\begin{lemma}
	\label{lemma:2log2}
	Let $k \in \N$.
	Then
	\begin{equation*}
		\frac{1}{k} \exp \paren[auto](){-\frac{2 \log(2)}{k + 1}}
		\le
		\frac{1}{k+1}
		,
	\end{equation*}
	where $\log$ denotes the natural logarithm.
\end{lemma}
\begin{proof}
	The claim is equivalent to:
	\begin{equation*}
		2 \log(2)
		\ge
		(k+1) \log \paren[auto](){1 + \frac{1}{k}}
		.
	\end{equation*}
	Consider the function $f(x) = (x + 1) \log\paren[auto](){1 + \frac{1}{x}}$ for $x > 0$ with derivative
	\begin{equation*}
		f'(x)
		=
		\log \paren[auto](){1 + \frac{1}{x}} - \frac{1}{x}
		\le
		0
		.
	\end{equation*}
	Thus, for $k \in \N$, we have $2 \log(2) = f(1) \ge f(k) = (k+1)\log\paren[auto](){1 + \frac{1}{k}}$.
\end{proof}

\begin{lemma}
	\label{lemma:harmonic-descent}
	Consider sequences $(d_k)$ and $(c_k)$ such that $d_k > 0$ and $c_k \ge c > 0$ hold for all $k$.
	Suppose moreover that for all $k \in \N$, we have
	\begin{equation*}
		d_{k+1}
		\le
		d_k \exp(-c_k \, d_k)
		.
	\end{equation*}
	Then there exists $k_0 \in \N$ such that for any $k > k_0$,
	\begin{equation*}
		d_k
		\le
		\frac{2 \log(2)}{c \, (k - k_0)}
		.
	\end{equation*}
\end{lemma}
\begin{proof}
	We provide a proof by induction.
	Suppose that $d_k \le \frac{2 \log (2)}{c \, (k - k_0)}$ for some $k > k_0$ for some $k_0 \in \N$.
	Then either $d_k \le \frac{2 \log(2)}{c \, (k+1 - k_0)}$, and we directly conclude $d_{k+1} \le d_k \underbrace{\exp\paren[auto](){-c_k \, d_k}}_{\le 1} \le \frac{2 \log(2)}{c \, (k+1 - k_0)}$, or else we have the bound
	\begin{equation*}
		d_{k+1}
		\le
		\frac{2 \log(2)}{c \, (k - k_0)} \exp\paren[auto](){-\frac{2 \log(2)}{k+1 - k_0}}
		\le
		\frac{2 \log(2)}{c \, (k+1 - k_0)}
	\end{equation*}
	due to \cref{lemma:2log2}.

	It remains to show that there exists an initial $k_0 \in \N$ such that $d_{k_0 + 1} \le \frac{2 \log(2)}{c}$ holds.
	Suppose, to the contrary, that $d_k > \frac{2 \log(2)}{c}$ holds for all $k \in \N$.
	Then
	\begin{equation*}
		d_{k+1}
		\le
		d_k \exp(-2 \log(2))
		=
		\frac{d_k}{4}
		.
	\end{equation*}
	By induction, we obtain $d_k \le \frac{d_0}{4^k}$.
	For $k \ge \frac{1}{\log(4)} \log \paren[big](){\frac{c \, d_0}{2 \log(2)}}$, this implies $d_k \le \frac{2 \log(2)}{c}$.
	Consequently, the claim holds with $k_0 \le \floor[big]{\frac{1}{\log(4)} \log \paren[big](){\frac{c \, d_0}{2 \log(2)}}}$.
\end{proof}

Our main result for this section is given in the following theorem.
We provide convergence results for SGD with the step sizes from \eqref{eq:perfect-step-size} in the interpolating and non-interpolating cases.
\begin{theorem}
	\label{theorem:SGD-global-convergence}
	Consider a $(\mu, L)$-feasible SOP, and denote $V_0 = \E[big]{\xi}{\norm{f_\xi'(w^\star)}_{X^*}^2}$ and
	\begin{equation*}
		C
		\coloneqq
		\sup_{k \in \N} \E[big]{\xi_k}{\norm{f_{\xi_k}'(w_k)}_{X^*}^2}
		.
	\end{equation*}
	Then, using SGD with step sizes determinded by \eqref{eq:perfect-step-size}, we obtain:
	\begin{enumerate}
		\item\label{item:harmonic-convergence}
			In case $V_0 > 0$, there exists $k_0 \in \N$ such that
			\begin{equation*}
				\E[auto]{0:k-1}{D_k}
				\le
				\frac{L}{\mu^2} \frac{C}{k - k_0}
			\end{equation*}
			for all $k > k_0$.

		\item\label{item:linear-convergence}
			In case $V_0 = 0$, and the SOP is pointwise $(\mu_\xi, L_\xi)$-feasible, then we have
			\begin{equation*}
				\E[auto]{0:k-1}{D_k}
				\le
				\theta^k D_0
			\end{equation*}
			with $\theta = 1 - \frac{\mu^2}{2 \, L \, \Lmax}$.
	\end{enumerate}
	\begin{proof}
		\begin{enumerate}
			\item
				From \eqref{eq:linear-convergence-mu-alpha} we infer:
				\begin{equation}
					\label{eq:bound-on-Dk:1}
					\E[auto]{0:k}{D_{k+1}}
					\le
					\E[auto]{0:k-1}{D_k (1 - \mu \, \alpha_k)}
					.
				\end{equation}
				For the step sizes from \eqref{eq:perfect-step-size} we conclude
				\begin{equation*}
					\alpha_k
					=
					\frac{\norm{F'(w_k)}_{X^*}^2}{L \, \E[big]{\xi_k}{\norm{f_{\xi_k}'(w_k)}_{X^*}^2}}
					\ge
					\frac{\norm{F'(w)}_{X^*}^2}{C L}
					\ge
					\frac{2 \, \mu D_k}{C L}
				\end{equation*}
				and thus by \eqref{eq:bound-on-Dk:1}
				\begin{equation*}
					\E[auto]{0:k}{D_{k+1}}
					\le
					\E[auto]{0:k-1}{D_k\paren[auto](){1 - \frac{2 \, \mu^2}{C L} \, D_k}}
					.
				\end{equation*}
				Now noticing that $x \mapsto x \, (1 - cx)$ is concave in $x$ for any non-negative~$c$, we have by Jensen's inequality (\cf \cite[Theorem~1.5.1]{Durrett:2019:1}):
				\begin{align*}
					\E[auto]{0:k-1}{D_k\paren[auto](){1 - \frac{2 \, \mu^2}{C L} \, D_k}}
					&
					\le
					\E[auto]{0:k-1}{D_k}\paren[auto](){1 - \frac{2 \, \mu^2}{C L} \E[auto]{0:k-1}{D_k}}
					\\
					&
					\le
					\E[auto]{0:k-1}{D_k} \exp\paren[auto](){-\frac{2 \, \mu^2}{C L} \E[auto]{0:k-1}{D_k}}
					.
				\end{align*}
				Thus, denoting $d_k \coloneqq \E[auto]{0:k-1}{D_k}$, we have shown:
				\begin{equation*}
					d_{k+1}
					\le
					d_k \exp\paren[auto](){-\frac{2 \, \mu^2}{C L} \, d_k}
					.
				\end{equation*}
				Now the results follows from \cref{lemma:harmonic-descent}.

			\item
				We have
				\begin{equation*}
					\E[auto]{\xi_k}{\norm{{f_{\xi_k}'(w_k)}}_{X^*}^2}
					\le
					4 \, \Lmax \, D_k + 2 \, V_0
				\end{equation*}
				due to \cref{lemma:adapted-variance-bound}. Thus, using $V_0 = 0$, $\norm{F'(w)}^2 \ge 2 \, \mu D_k$ we obtain:
				\begin{equation*}
					\alpha_k
					=
					\frac{\norm{F'(w_k)}_{X^*}^2}{L \, \E[big]{\xi_k}{\norm{{f_{\xi_k}'(w_k)}}_{X^*}^2}}
					\ge
					\frac{2 \, \mu D_k}{4 \, L \, \Lmax \, D_k}
					\ge
					\frac{\mu}{2 \, L \, \Lmax}
					.
				\end{equation*}
				Again, denoting $d_k = \E[auto]{0:k-1}{D_k}$, we obtain from \eqref{eq:linear-convergence-mu-alpha}
				\begin{equation*}
					d_{k+1}
					\le
					\E[big]{0:k-1}{(1 - \mu \, \alpha_k) \, D_k}
					\le
					\paren[auto](){1 - \frac{\mu^2}{2 \, L \, \Lmax}} d_k
				\end{equation*}
				and thus the claimed linear convergence.
		\end{enumerate}
	\end{proof}
\end{theorem}

\begin{remark}
	In the setting of \cref{theorem:SGD-global-convergence} the step size $\alpha_k$ is also a random variable, which depends on $\xi_0, \dots, \xi_{k-1}$.
	This complicates the proof for the first case, and motivates the incorporation of the assumption $\sup_{k \in \N}{ \E[big]{\xi_k}{\norm{f_{\xi_k}'(w_k)}_{X^*}^2}} < \infty$, as well as the usage of Jensen's inequality.
	The assumption is common in the literature dealing with convergence results for SGD.
	In any practical scenario an assumption like this seems inevitable (\cite{NguyenNguyenVanDijkRichtarikScheinbergTakac:2018:2,JiangStich:2023:1}).
\end{remark}
\begin{remark}
	The linear convergence established in \cref{theorem:SGD-global-convergence}~\ref{item:linear-convergence} is based on a lower bound on the step size used, which in turn is a consequence of an upper bound on the variance obtained from \cref{lemma:adapted-variance-bound}.
	Thus, the pointwise $(\mu_\xi, L_\xi)$-feasibility assumption, in particular the convexity of almost every $f_\xi$, can be dropped in favor of weaker assumptions that still yield bounds on the variance.
	One possible assumption is the setting of \cref{lemma:first-bound-on-variance}.
	This, of course, leads to a larger value of~$\theta$.
\end{remark}

\section{Estimation Techniques}
\label{section:estimation-techniques}

This section provides an outline of the estimation methods used to obtain the quantities required by our adaptive schemes.
In particular, we intend to explain the ideas behind \cref{algorithm:rough-sketch} stated below.
The proposed estimators described in this section require an additional piece of information that is not available during a traditional SGD run.
This extra information is the additional evaluation of the $k$-th sampled function (without gradient) at the next iterate $w_{k+1}$, \ie, we additionally evaluate
\begin{equation*}
	f_{\xi_k}(w_{k+1})
	,
\end{equation*}
after performing the step $w_{k + 1} = w_k - \alpha_k\nabla f_{\xi_k}(w_k)$ in the $k$-th iteration.
In classical machine learning settings, this comes at the moderate cost of one additional forward pass per iteration, which roughly increases the cost per iteration about \SI{33}{\percent}.

Based on the results of \cref{section:convergence-to-optimality}, the following SGD algorithm with adaptive step size control is proposed.
A more detailed version is provided in \cref{section:complete-algorithm-in-pseudo-code}.
\begin{algorithm}[htp]
	\caption{SGD with adaptive step size control}
	\label{algorithm:rough-sketch}
	\begin{algorithmic}[1]
		\State Find some $\alpha_0 > 0$ via line search such that $f_{\xi_0}(w_0 - \alpha \nabla f_{\xi_0}(w_0)) < f_{\xi_0}(w_0)$
		\State Initialize $L_0 \gets \frac{1}{\alpha_0}$
		\State Initialize $\sigma_0^2 \gets 0$ and $g_0 \gets \norm{f_{\xi_0}'(w_0)}_{X^*}^2$
		\For{$k \ge 1$}
		\State $w_{k+1} \gets w_k - \alpha_k \nabla f_{\xi_k}(w_k)$
		\Comment{classical SGD step}
		\State Evaluate $f_{\xi_k}(w_{k+1}$)
		\Comment{additional information}
		\If{$f_{\xi_k}(w_{k+1}) > f_{\xi_k}(w_k)$}
		\State $w_{k+1} \gets w_k$
		\Comment{reject the step}
		\State Continue with next $k$, perform line search on the next sample $\xi_{k+1}$
		\EndIf
		\State Update the estimators for $L_{k+1}$, $\sigma_{k+1}^2$ and $g_{k+1}$, as described below
		\State Select the next step size $\alpha_{k+1}$ according to \eqref{eq:perfect-step-size}, using estimators $L$, $\sigma^2$, $g$
		\EndFor
	\end{algorithmic}
\end{algorithm}
In this algorithm, $L_k$, $\sigma_k^2$ and $g_k$ represent the estimators for the nonlinearity~$L$, the local variance $\Var{\xi_k}{f_\xi'(w)}$, and $\E{\xi}{\norm{f_\xi'(w)}_{X^*}^2}$, respectively.

The remainder of this section is organized as follows.
We first introduce the method of \emph{exponential smoothing}, which is a widely used technique in signal processing to obtain a moving average over a time series.
This method is employed in our work to average the individual observations of our estimators in order to obtain a usable, relatively stable estimator for each quantity.
We then describe how to compute the individual observations of our estimators in each iteration.

\subsection{Exponential Smoothing}
\label{subseciton:exponential-smoothing}

Exponential smoothing techniques are well known in the signal processing literature.
An unknown, potentially time-dependent and noisy quantity is estimated by previous observations and one new observation.
Let $\tilde x_k$ be the $k$-th observation of the quantity of interest.
The smoothed estimate $x_k$ is then given by
\begin{equation*}
	x_k
	=
	\gamma \, x_{k-1} + (1 - \gamma) \, \tilde x_k
\end{equation*}
for some discount factor $\gamma \in (0,1)$.

A larger value of $\gamma$ reduces the impact of the most recent observation, resulting in a smoother estimate, but also increases the delay in the estimation.
By contrast, using a smaller value of $\gamma$ results in a more responsive yet more noisy estimation.

A simple extension of this approach is to make $\gamma$ time-dependent, for example one could use $ x_k = \gamma_k \, x_{k-1} + (1 - \gamma_k) \, \tilde x_k$ with $\gamma_k = 1 - \frac{1}{k}$.
This approach yields the classical average $x_k = \frac{1}{k} \sum_{i = 1}^k \tilde x_i$.
We employ time-dependent discount factors to smooth our estimators, using $\gamma_k = 1 - k^{-\eta}$ for some $\eta \in [\frac{1}{2} ,1)$ for all estimators.\footnote{In our numerical experiments we used $\eta = 0.7$ for quadratic problems (\cref{subsection:quadratic-SOPs}) and $\eta = 0.8$ for image classification (\cref{subsection:image-classification}). We observed that the performance of the algorithm is not sensitive to the selection of $\eta$.}

This results in an estimation approach that is able to quickly adapt to the problem in the early stages of the algorithm while providing a more stable estimate in later stages.
Stable estimates are observed to be crucial to obtain convergence.
This can be explained by the representation of the step size in equation \eqref{eq:perfect-step-size}, where the difference between the two estimated values $\E{\xi_k}{\norm{f_\xi'(w_k)}_{X^*}^2}$ and $\Var{\xi_k}{f_\xi'(w_k)}$ plays a central role.
If this difference does not tend to zero, the step size will not tend to zero.
Thus, decreasing the noise in the estimates is crucial.
Fortunately, this can be achieved by adjusting $\gamma_k$ without any additional computational complexity.

Exponential smoothing techniques are not new to the stochastic optimization literature.
For instance, the averaging over the second moments of the gradient found in Adam (\cite{KingmaBa:2015:1}) can be regarded as an exponential smoothing technique.
Additionally, classical momentum schemes can also be considered as exponential smoothing approaches for the gradient if the hyperparameters are properly chosen.

\subsection{Nonlinearity Estimation}
\label{subsection:nonlinearity-estimation}

\subsubsection{Description of the Estimator}
\label{subsubsection:description-of-the-estimator}

One of the quantities of large interest in step size selection is the nonlinearity of the problem, expressed through the Lipschitz constant~$L$ of the gradient.
In this section, we provide a method to estimate this quantity.
We further discuss the computational complexity involved.
It is well known that $L$ is a Lipschitz constant for the gradient of~$F$ if and only if
\begin{equation*}
	F(w + \delta_w)
	\le
	F(w) + F'(w) \, \delta_w + \frac{L}{2} \norm{\delta_w}_X^2
	\quad
	\text{for all }
	w, \delta_w \in X
	.
\end{equation*}
Rearranging the terms gives
\begin{equation*}
	L
	\ge
	2 \, \frac{F(w + \delta_w) - F(w) - F'(w) \, \delta_w}{\norm{\delta_w}^2}
	.
\end{equation*}
Note that if $\delta_w = -\alpha \nabla F(w)$, we have $F'(w) \, \delta_w = -\alpha \, \norm{\nabla F(w)}^2$ and $\norm{\delta_w}^2 = \alpha^2 \, \norm{\nabla F(w)}^2$.
Motivated by this, we evaluate
\begin{equation*}
	\tilde L_k
	=
	2 \, \frac{f_{\xi_k}(w_{k+1}) - f_{\xi_k}(w_k) + \alpha_k \, \norm{\nabla f_{\xi_k}(w_k)}_X^2}{\alpha_k^2 \, \norm{\nabla f_{\xi_k}(w_k)}_X^2}
\end{equation*}
in the $k$-th iteration.
We use this as an individual observation and then apply exponential smoothing on these observations.
That is, in each iteration we correct the estimate~$L_k$ for the nonlinearity as follows,
\begin{equation*}
	L_{k+1}
	\coloneqq
	\gamma_k \, L_k + (1 - \gamma_k) \, \tilde L_k
	.
\end{equation*}
We thus only need an initial value~$L_0$.
In \cref{subsection:practical-aspects:initialization} we describe a way to select $L_0$.

\subsubsection{Discussion of the Estimator}

It is clear that the method presented here underestimates the theoretical value of~$L$.
Rather, our method focuses on an average over all Lipschitz values for $\nabla F(w)$ in the directions SGD takes.
In practice this is adequate for the decrease condition to hold and, therefore, enough for obtaining descent.
Clearly, this argument leans more towards a heuristic approach.

A more conservative estimate of the problem's nonlinearity would be to use the maximum of all estimates of~$L$.
However, in practice this results in exceedingly large values for $L_k$, since it also includes outliers and the method does not adjust well to areas with smaller Lipschitz constants.
Our numerical experiments have shown that the latter method produces a step size that is overly pessimistic, since it employs a global estimate for a quantity that only acts locally.

\subsubsection{Computational Cost}
\label{subsubsection:computational-cost}

Besides the evaluation of $f_{\xi_k}(w_{k+1})$, additional computational costs arise from the evaluation of $\norm{\nabla f_{\xi_k}(w_k)}_X^2$.
Usually, $f'_{\xi_k}(w_k)$ and $\nabla f_{\xi_k}(w_k)$ are computed during the (preconditioned) SGD step.
In classical machine learning tasks, $f_{\xi_k}'(w_k)$ is computed via back-propagation, then possibly a preconditioner is applied to obtain $\nabla f_{\xi_k}(w_k) = \inv{H} f_{\xi_k}'(w_k)$.
It follows that
\begin{equation*}
	\norm{\nabla f_{\xi_k}(w_k)}_X^2
	=
	\dual{f_{\xi_k}'(w_k)}{\nabla f_{\xi_k}(w_k)}
	.
\end{equation*}
Thus, evaluation of $\norm{\nabla f_{\xi_k}(w_k)}$ comes at the cost of one duality product, \ie, one inner product of the vectors $f_{\xi_k}'(w_k)^\transp$ and the gradient $\nabla f_{\xi_k}(w_k)$.

\subsection{Variance Estimation}
\label{subsection:variance-estimation}

Consider an SGD step performed with step size $\alpha_k$.
Note that when $\alpha_k$ is sufficiently small, we have
\begin{equation*}
	\E{\xi_k}{f_{\xi_k}(w_{k+1})}
	<
	\E{\xi_{k+1}}{f_{\xi_{k+1}}(w_k)}
	=
	F(w_{k+1})
	,
\end{equation*}
which provides a biased estimate of the true functional value at the iterate~$w_{k+1}$.
This is because the search direction $\delta_k$ is selected to minimize $f_{\xi_k}$, not $F$.
However, by comparing the unbiased estimator $f_{\xi_{k+1}}(w_{k+1})$ to the biased estimator $f_{\xi_k}(w_{k+1})$, we can determine a notion of the local variance.

In order to quantify the above heuristic, recall that for sufficiently smooth functions we have
\begin{equation*}
	f(w + \delta_w)
	=
	f(w) + f'(w)\delta_w + O(\norm{\delta_w}^2)
	.
\end{equation*}
Applying this expansion to $f_{\xi_k}(w_{k+1})$ and $f_{\xi_{k+1}}(w_{k+1})$, we obtain
\begin{align*}
	f_{\xi_k}(w_{k+1})
	&
	=
	f_{\xi_k}(w_k) - \alpha_k \, \norm{f_{\xi_k}'(w_k)}_{X^*}^2 + O(\alpha_k^2)
	\intertext{and}
	f_{\xi_{k+1}}(w_{k+1})
	&
	=
	f_{\xi_{k+1}}(w_k) - \alpha_k \, \dual{f_{\xi_{k+1}}'(w_k)}{\nabla f_{\xi_k}(w_k)} + O(\alpha_k^2)
	.
\end{align*}
We thus obtain
\begin{multline}
	f_{\xi_{k+1}}(w_{k+1}) - f_{\xi_k}(w_{k+1})
	\\
	=
	f_{\xi_{k+1}}(w_k) - f_{\xi_k}(w_k)
	+
	\alpha_k \, \norm{f_{\xi_k}'(w_k)}_{X^*}^2
	-
	\alpha_k f_{\xi_{k+1}}'(w_k)\nabla f_{\xi_k}(w_k)
	+
	O(\alpha_k^2)
	\label{eq:difference explanation variance}
\end{multline}
for the difference.
On the one hand, we have $\E{\xi_{k+1}}{f_{\xi_{k+1}}(w_k)} = \E{\xi_k}{f_{\xi_k}(w_k)}$.
On the other hand, since $\xi_k$ and $\xi_{k+1}$ are independent, we obtain
\begin{equation*}
	\E{\xi_k, \xi_{k+1}}{\dual{f_{\xi_{k+1}}'(w_k)}{\nabla f_{\xi_k}(w_k)}}
	=
	\E{\xi_{k+1}}{f_{\xi_{k+1}}'(w_k)} \, \E{\xi_k}{\nabla f_{\xi_k}(w_k)}
	=
	\norm{F'(w_k)}_{X^*}^2
	.
\end{equation*}
We thus get, taking the expectation of \eqref{eq:difference explanation variance}:
\begin{align*}
	\E{\xi_k, \xi_{k+1}}{f_{\xi_{k+1}}(w_{k+1}) - f_{\xi_k}(w_{k+1})}
	&
	=
	\alpha_k \, \E{\xi_k}{\norm{f_{\xi_k}'(w_k)}^2}_{X^*}
	-
	\alpha_k \, \norm{F'(w_k)}_{X^*}^2 + O(\alpha_k^2)
	\\
	&
	=
	\alpha_k \Var{\xi_k}{f_{\xi_k}'(w_k)} + O(\alpha_k^2)
	.
\end{align*}
Neglecting the second order term, we obtain a way to estimate the variance
\begin{equation*}
	\E{\xi_k, \xi_{k+1}}{f_{\xi_{k+1}}(w_{k+1}) - f_{\xi_k}(w_{k+1})}
	\approx
	\alpha_k \Var{\xi_k}{f_{\xi_k}'(w_k)}
	.
\end{equation*}
Motivated by this, in the $(k+1)$-th iteration, after evaluating $f_{\xi_{k+1}}(w_{k+1})$, we evaluate
\begin{equation*}
	\tilde \sigma_k^2
	=
	\frac{f_{\xi_{k+1}}(w_{k+1}) - f_{\xi_k}(w_{k+1})}{\alpha_k}
\end{equation*}
and then again use exponential smoothing to update our variance estimation
\begin{equation*}
	\sigma_{k+1}^2
	=
	\gamma_k \, \sigma_k^2 + (1 - \gamma_k) \, \tilde \sigma_k^2
	.
\end{equation*}
Again, this scheme needs an initialization, we comment on this in \cref{subsection:practical-aspects:initialization}.

\subsection{\texorpdfstring{Accessibility of $\E[big]{\xi_k}{\norm{f_{\xi_k}'(w_k)}_{X^*}^2}$}{Accessibility of the Expected Value of the Gradient Norm}}
\label{subsection:accessibility-of-expected-value-gradient-norm}

The remaining quantity in \eqref{eq:perfect-step-size} to be estimated is
\begin{equation*}
	\E{\xi_k}{\norm{f_{\xi_k}'(w_k)}_{X^*}^2}
	.
\end{equation*}
We also use exponential smoothing to estimate this quantity.
To this end, in every iteration we evaluate
\begin{equation*}
	\tilde g_k
	=
	\norm{f_{\xi_k}'(w_k)}_{X^*}^2
\end{equation*}
and then update
\begin{equation*}
	g_{k+1}
	=
	\gamma_k \, g_k + (1 - \gamma_k) \, \tilde g_k
	.
\end{equation*}
This computation has negligible cost, since $\tilde g_k = \norm{f_{\xi_k}'(w_k)}_{X^*}^2$ has already been evaluated to obtain $\tilde L_k$.
Again, for initialization we refer to \cref{subsection:practical-aspects:initialization}.

\section{Practical Aspects}
\label{section:practical-aspects}

In this section, we provide some remarks on aspects that may affect the performance of our method in practical scenarios.
This is particularly important when solving problems beyond our theoretical framework.
Specifically, we concentrate on initialization (which applies to problems within and beyond our theory), as well as on safeguards against unreliable estimates and global convergence for nonconvex problems.

\subsection{Initialization}
\label{subsection:practical-aspects:initialization}

Obviously, all our estimation schemes need an initial value.
We use the following strategies:
\begin{itemize}
	\item
		For $\alpha$ and $L$, we run a line search on the very first training sample $f_{\xi_0}$, looking for a step size $\alpha_0 > 0$ satisfying $f_{\xi_0}(w_1) \le f_{\xi_0}(w_0)$ and then using $L_0 = \frac{1}{\alpha_0}$.
		This allows users to start the algorithm without any knowledge about the scaling of the problem.
		In the preceding iterates, due to relatively small values of $\gamma$, the algorithm is able to quickly adapt.

	\item
		We use $0$ as initialization for the variance.
		This is motivated by the fact that, from the results in \cref{subsection:variance-bounds-and-initial-convergence}, we do not need to take the true local variance into account at the beginning of the algorithm.

	\item
		We use $\norm{f_{\xi_0}'(w_0)}_{X^*}^2$ as initialization for the $\E[big]{\xi}{\norm{f_\xi'(w_k)}_{X^*}^2}$ estimate.
\end{itemize}

\subsection{Safeguards}

All estimated quantities are subject to noise, so besides exponential averaging, some safeguards are in order.
Some glitches are easy to spot and avoid.
For example, it holds
\begin{equation*}
	\E[big]{\xi}{\norm{f_\xi'(w)}_{X^*}^2}
	=
	\norm{F'(w)}_{X^*}^2 + \Var{\xi}{f_\xi'(w)}
\end{equation*}
and therefore
$\E[big]{\xi}{\norm{f_\xi'(w)}_{X^*}^2} - \Var{\xi}{f_\xi'(w)} \ge 0$.
In case the corresponding estimated quantities violate this inequality, the suggested learning rate from the corresponding iteration can be disregarded, since it would be negative anyway.

Also, since $\Var{\xi}{f_\xi'(w)} \ge 0$, we can ignore negative values of the corresponding estimates.

Further, we could impose upper and lower limits for all estimates in order to avoid exploding values in settings that are outside the realm of our theory.
Furthermore, in case we observe $f_{\xi_k}(w_{k+1}) > f_{\xi_k}(w_k)$, we reject the step and proceed with $w_{k+1} = w_k$.
Since this observation is likely due to a too large step size, we perform a line search \emph{in the next iteration} in this case.
In our numerical experiments we observed this to happen only occasionally during runs, mostly at the beginning of the run of our algorithm, where the estimates are not reliable yet.

\subsection{Global and Local Phases for Nonconvex Problems}
\label{subsection:global-local-phases}

As classical machine learning problems are highly nonconvex, we start our algorithm with a global phase, where we exclude the local variance from the step size calculation, to prevent the step size from becoming too small too early.
This approach is essential to move away from the initialization point and avoid getting trapped in local minimizers.
In this way we are able to leverage this well known property of SGD.
After a certain number of iterations, we enter the \emph{local phase}, during which we additionally consider the variance.
Intuitively, this means that we have reached a neighborhood of the solution and start to fine tune the iterates.

\section{Numerical Results}
\label{section:numerical-results}

To illustrate the theoretical results given above, numerical experiments were performed, where we used our algorithm outlined in \cref{algorithm:rough-sketch} to solve several optimization problems.
In the first experiment, we focus on a problem setting that meets the assumption of our theory.
Here, we construct SOPs from quadratic functions.
The resulting SOPs are $(\mu, L)$-feasible, and for sufficiently small noise level $\sigma_A$, see below, also pointwise $(\mu_\xi, L_\xi)$-feasible.
In the second experiment, we consider four different image classification tasks with ReLU networks.
Obviously, our assumptions are not met here, as ReLU networks are not differentiable, and the corresponding target function is known to be highly nonconvex.
Remarkably, our method is still able to adapt to the underlying problem and provides good step sizes.

No preconditioning is used in any of the experiments.
Consequently, we measure convergence in the Euclidean norm $\norm{\cdot}_2$ of $\R^n$.

\subsection{Quadratic SOPs}
\label{subsection:quadratic-SOPs}

Given an orthogonal matrix $S \in \R^{n\times n}$ and a diagonal matrix $D = \diag(\lambda_1, \dots, \lambda_n)$, we construct an SOP as follows.
We set the mean Hessian to $A \coloneqq S^\transp D S$ and select a noise level $\sigma_A > 0$.
In every iteration, we sample a random matrix $\Xi \in \R^{n \times n}$ with every entry $\xi_{ij}$ drawn from the uniform distribution on $[-\sigma_A,\sigma_A]$.
Then we let $W \coloneqq \Xi^\transp \Xi - \frac{2}{3} \sigma_A^3 \id$.
As is easily checked, this ensures $\E{\Xi}{W} = 0$.
We then use $A_\xi = A + W$ as the matrix for the quadratic SOP in the respective iteration.

For $b \in \R^n$, we choose a noise level $\sigma_b \ge 0$ and sample every entry of $b_\xi$ from the uniform distribution on $[-\sigma_b, \sigma_b]$.

We then consider the problem to minimize the expected value of
\begin{equation}
	\label{eq:quadratic-sop:objective}
	f_\xi(w)
	=
	\frac{1}{2} w^\transp A_\xi w
	+
	b_\xi^\transp w
	.
\end{equation}
In each scenario (interpolating and non-interpolating), we perform ten different test runs using different random seeds and plot the average and one standard deviation (of the $\log_{10}$ of the respective quantity) in our plots.

\subsubsection{Non-Interpolating Case}
\label{subsubsection:non-interpolating-case}

The non-interpolating case corresponds to the case $\E[big]{\xi}{\norm{\nabla f_\xi(w^\star)}_2^2} > 0$ and thus to $\sigma_b > 0$.
In this situation, we expect the step size to descend to zero and convergence of $\norm{w_k - w^\star}_2$ to zero at a rate of at least $O(1/\sqrt k)$ due to \cref{theorem:SGD-global-convergence}.

We tested two different possible variants for~$D$.
We choose $L$ and $\mu$ and eigenvalues $\lambda_i = \mu + (\frac{i-1}{n-1})^2 \, (L - \mu)$ for $1 \le i \le n$.
In the first scenario, we fix $\mu = 1$ and let $L$~grow.
In the second scenario, we fix $L = 1$ and let $\mu$~decrease.
In both scenarios, the condition number $L/\mu$ varies between $10$ and $10^4$.
We choose the problem dimension to be $n = 50$.

\begin{figure}[ht]
	\centering
	\hfill
	\includegraphics[width=\textwidth]{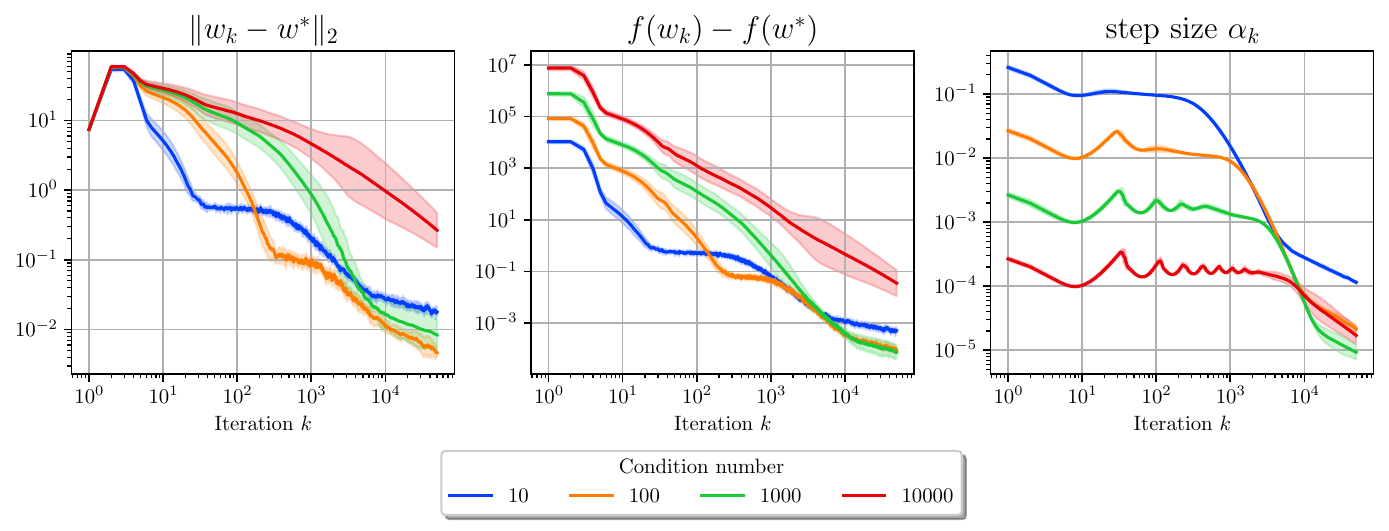}
	\caption{Non-interpolating case: performance of adaptive step size control for the first scenario ($\mu = 1$ and $L$~variable).}
	\label{figure:noninterpolating:variable-L}
\end{figure}

The results for the first scenario are shown in \cref{figure:noninterpolating:variable-L}.
We observe the expected behavior, \ie, step sizes which are proportional to~$L$ initially and then decrease to zero, which allows for convergence to the minimizer.

\begin{figure}[ht]
	\centering
	\hfill
	\includegraphics[width=\textwidth]{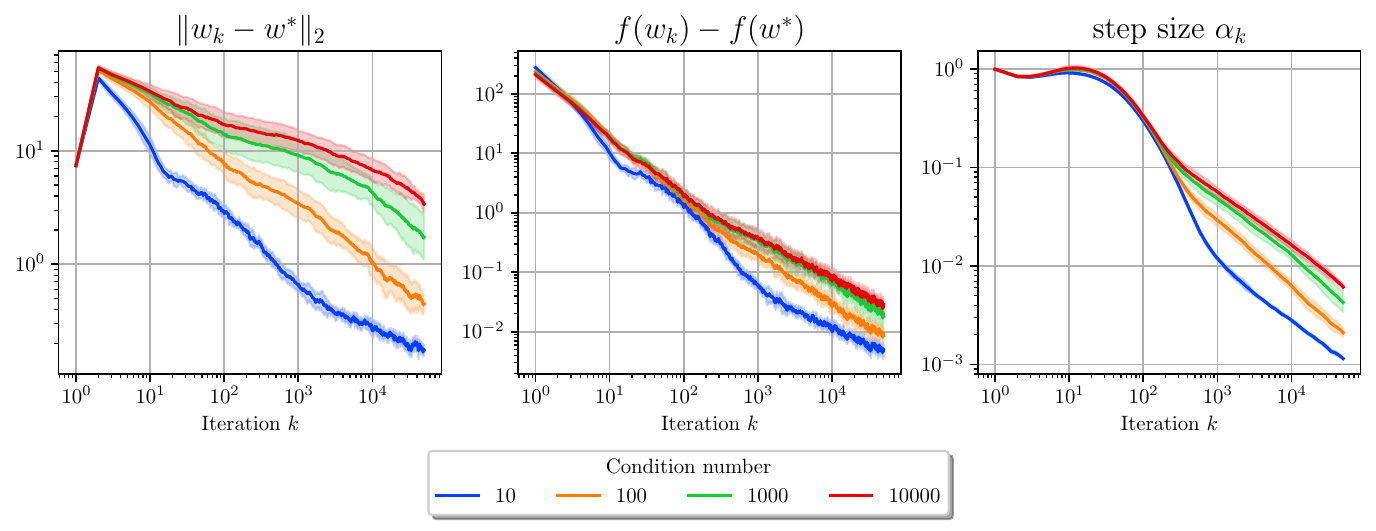}
	\caption{Non-interpolating case: performance of adaptive step size control for the second scenario ($L = 1$ and $\mu$~variable).}
	\label{figure:noninterpolating:variable-mu}
\end{figure}

The results for the second scenario are shown in \cref{figure:noninterpolating:variable-mu}.
Again, we observe the expected behavior, \ie, an initial step size proportional to $L = 1$, which then decreases.
This allows for convergence of the iterates and functional value of the target function.
Here, we observe that the convergence of the iterates depends more strongly on~$\mu$.
This is due to the small eigenvalues, which result in small gradients in the direction of the respective eigenspaces.

\subsubsection{Interpolating Case}
\label{subsection:interpolating-case}

The interpolating case corresponds to the absence of noise at the minimizer and thus to $\sigma_b = 0$.
In this situation, we expect linear convergence of the iterates and the functional value, and step sizes bounded from below.
Again, this is precisely the behavior we observe.
The oscillation in the step sizes can be explained by the estimation technique used to assess the nonlinearity ($L$).
It can easily be seen that our estimator for~$L$ always yields an estimate smaller than the true value of~$L$.
Thus, at some point, the step size might become too large, which leads to $f_{\xi_k}(w_{k+1}) > f_{\xi_k}(w_k)$ for some~$k$.
In this event, we reject the corresponding step and perform a line search on the next sample, which leads to a decreasing step size.
In later iterations, we use a larger averaging parameter (see \cref{subseciton:exponential-smoothing}), and thus damp the effect of too small estimates of~$L$.

\begin{figure}[ht]
	\centering
	\hfill
	\includegraphics[width=\textwidth]{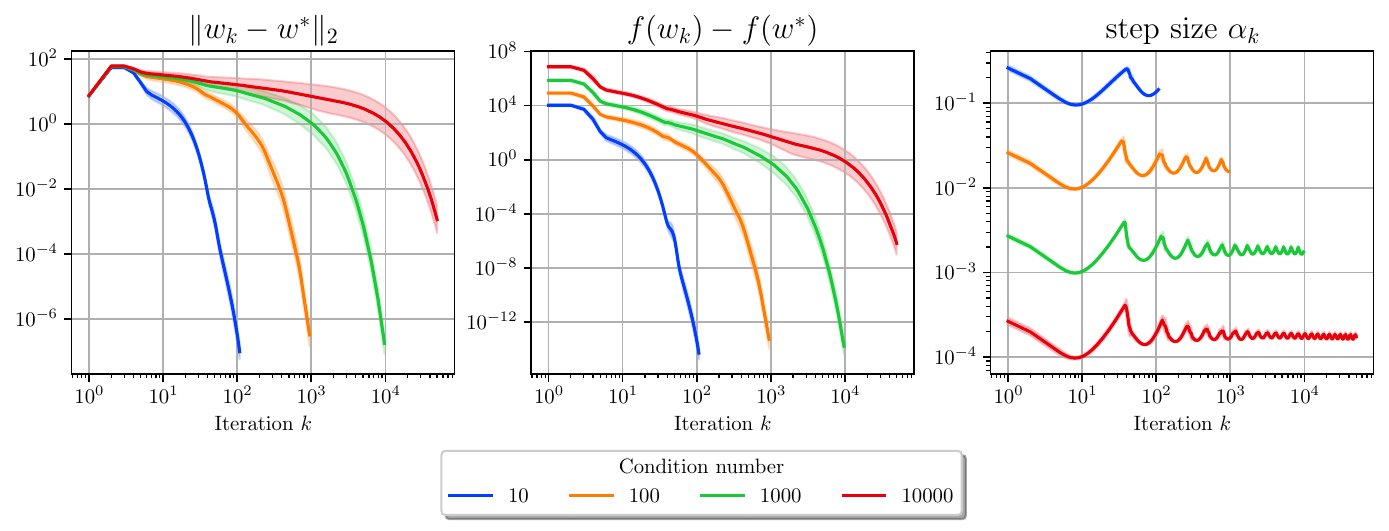}
	\caption{Interpolating case: performance of adaptive step size control for the first scenario ($\mu = 1$ and $L$~variable.)}
	\label{figure:interpolating:variable-L}
\end{figure}

\begin{figure}[H]
	\centering
	\hfill
	\includegraphics[width=\textwidth]{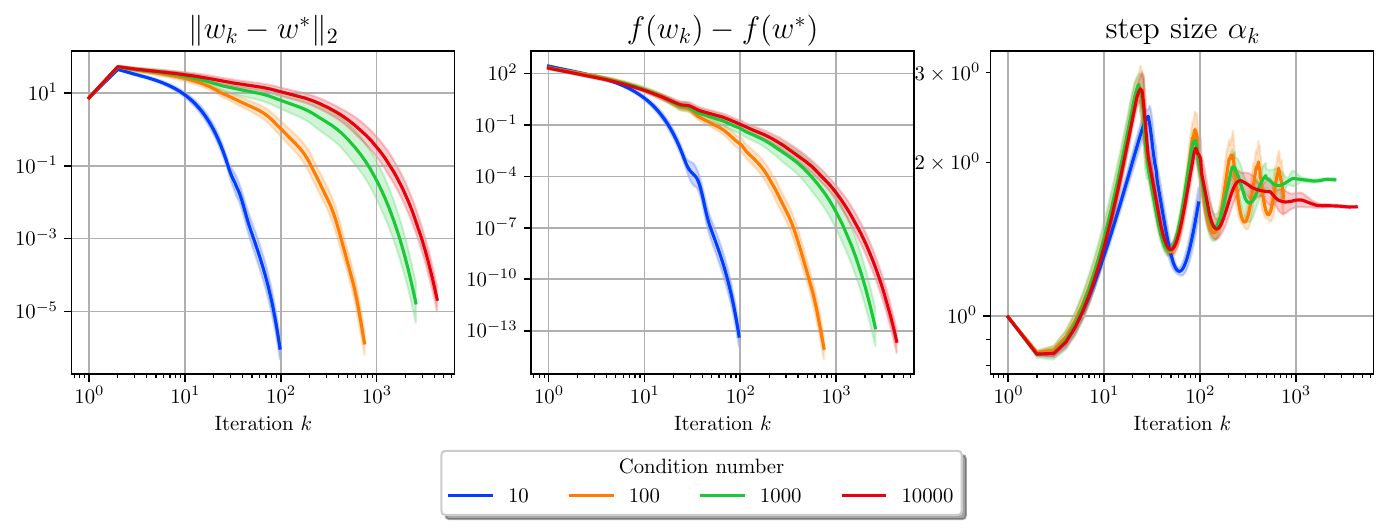}
	\caption{Interpolating case: performance of adaptive step size control for the second scenario ($L = 1$ and $\mu$~variable).}
	\label{figure:interpolating:variable-mu}
\end{figure}

\subsection{Image Classification Tasks}
\label{subsection:image-classification}

We also tested \cref{algorithm:rough-sketch} on four standard image classification tasks identified in the benchmarking paper \cite{SchmidtSchneiderHennig:2021:1}.
We used the benchmarking tool DeepOBS from \cite{SchneiderBallesHennig:2019:1} to test the performance of the algorithm.

The experiments were performed on a workstation equipped with an Intel i5-12500 CPU and NVIDIA RTX3070 GPU, using the \pytorch package.
A comprehensive overview of the problems we consider is given in \cref{table:test-problems}.
We refer the interested reader to \cite{SchneiderBallesHennig:2019:1} for details regarding the respective network architectures and to the original works \cite{XiaoRasulVollgraf:2017:1} (Fashion-MNIST), \cite{Krizhevsky:2009:1} (CIFAR-10 and CIFAR-100), \cite{NetzerWangCoatesBissaccoWuNg:2011:1} (SVHN) for the data sets.
The last column shows the approximate run time for one single training run.

Since all the problems considered here are highly nonconvex, we use the global-phase/local-phase approach discussed in \cref{subsection:global-local-phases}.
We start the local phase after \SI{60}{\percent} of the available epochs.

\begin{table}[ht]
	\begin{tabular*}{\textwidth}{@{\extracolsep{\fill}}llrrr}
		\multicolumn{5}{c}{Image Classification Tasks}
		\\
		\toprule
		Data set       & Model architecture      & Batch size & \#Epochs & Run time
		\\
		\midrule
		Fashion-MNIST  & Simple CNN: \emph{2c2d} & 128        & 100      & \SI{7}{\minute}
		\\
		CIFAR-10       & Simple CNN: \emph{3c3d} & 128        & 100      & \SI{25}{\minute}
		\\
		SVHN           & \emph{Wide ResNet 164}  & 128        & 160      & \SI{125}{\minute}
		\\
		CIFAR-100      & \emph{All-CNN-C}        & 256        & 350      & \SI{165}{\minute}
		\\
		\bottomrule
	\end{tabular*}
	\caption{Test problems considered for the image classification tasks.}
	\label{table:test-problems}
\end{table}

\begin{figure}[ht]
	\centering
	\hfill
	\includegraphics[width = \textwidth]{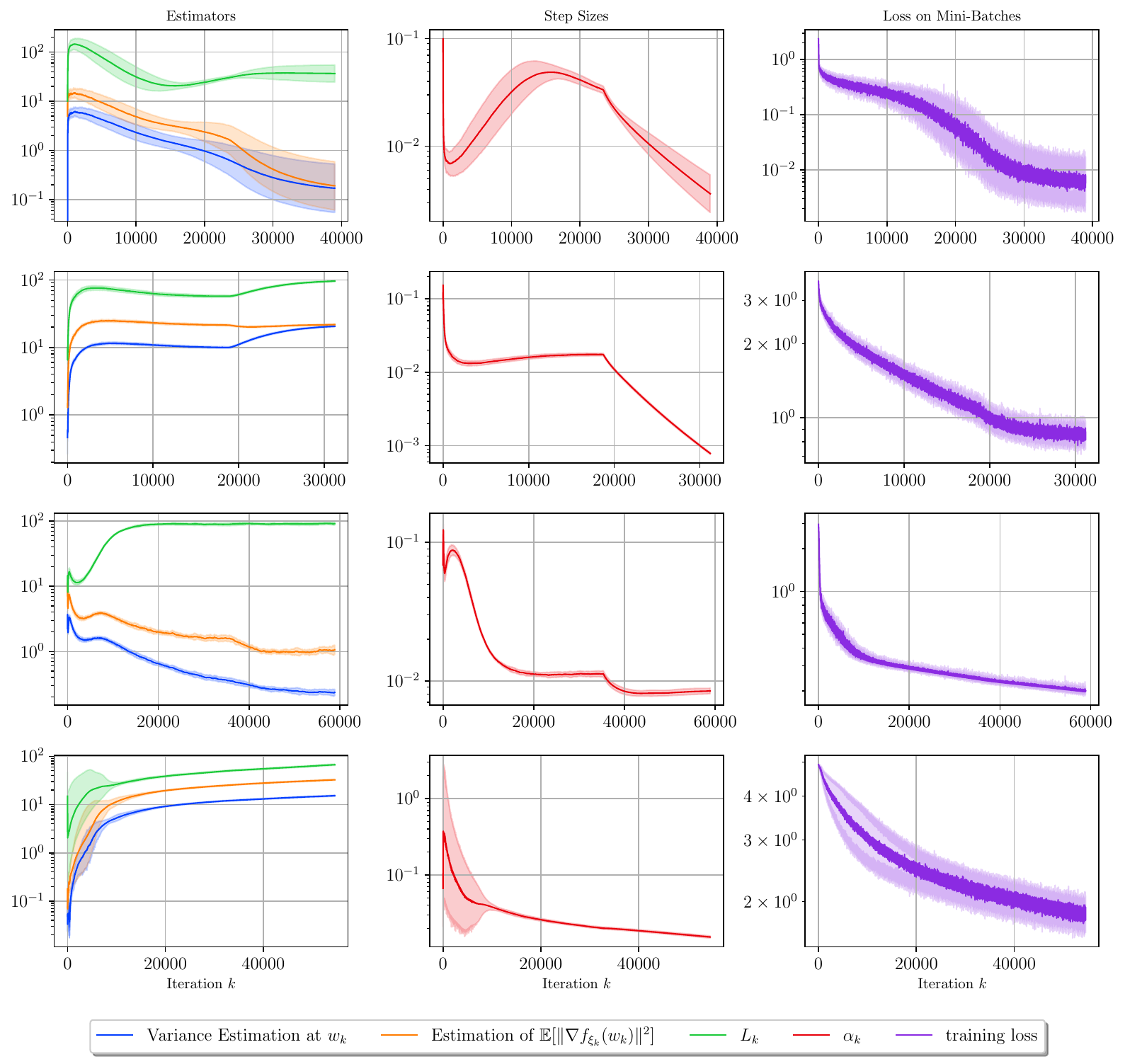}
	\caption{Performance on Image Classification data sets.
	Top row: Fashion-MNIST, second row: CIFAR-10, third row: SVHN, last row: CIFAR-100.}
	\label{figure:image-classification-data-sets}
\end{figure}

The results of our experiments are presented in \cref{figure:image-classification-data-sets}.
We use ten different random seeds for initialization and batch selection and plot the mean (solid line) as well as one standard deviation of the $\log_{10}$ of the respective value (shaded area).
\Cref{algorithm:rough-sketch} is clearly capable of adapting to these different settings.
The qualitatively different behavior of our estimators for the Lipschitz constant~$L$, the averaged norm of the gradient $\E{\xi_k}{\norm{\nabla f_{\xi_k}(w_k)}_X^2} = \E{\xi_k}{\norm{f_{\xi_k}'(w_k)}_{X^*}^2}$ and the variance $\Var{\xi_k}{f_\xi'(w)}$ across the different test problems demonstrates that the method is truly adaptive and provides problem specific step sizes for SGD.

\section{Discussion}
\label{section:discussion}

In this paper we introduced a novel technique to control the step sizes of SGD.
Our approach relies on the identification of computable quantities that we estimate during the run of the algorithm to obtain adaptive step sizes.
Our numerical experiments clearly show the adaptivity of our method.
We stress that (besides the global/local phase and the slightly different smoothing constants, see \cref{subseciton:exponential-smoothing}), no adjustment whatsoever to the algorithm was necessary for it to tackle all problems we considered --- quadratic SOPs as well as image classification tasks --- and no tuning of hyperparameters was necessary.
These advantages more than compensate the slightly increased cost for each iteration (see \cref{section:estimation-techniques}).

The numerical results are substantiated with a rigorous theoretical analysis of our scheme.
Under the assumption that the estimators we evaluate are precise, \cref{theorem:SGD-global-convergence} shows global convergence of SGD with the best known rates.

As a general observation, our algorithm finds step sizes that are comparable to good empirical choices found by a systematic search over several test runs.
This aligns with the main motivation of our research: to devise automatic step size adaptation in order to reduce the amount of manual tuning and repeated optimization runs.
Quantifying this observation by theoretical analysis and practical experiments is a current aspect of our research.

Up to now, our theoretical understanding is limited to strongly convex functions.
However, some results of our theory might carry over to a more general problem class.
Even in the absence of theory, our methods have proved to work well in practical experiments.

Our theoretical analysis clearly distinguishes between derivatives and gradients.
While in our numerical results we have only used the classical SGD method as a basic algorithm, our theoretical considerations include the option to use preconditioners, which will result in improved constants in the convergence results.
Our step size scheme, discussed in \cref{section:convergence-to-optimality}, is compatible with known preconditioning or second order methods that are used in the context of neural network training, specifically approximations to natural gradient methods or generalized Gauss-Newton methods as discussed for example in \cite{MartensSutskever:2012:1,Martens:2014:1,MartensGrosse:2015:1}.
Also the more general framework of \cite{HerzogKoehneKreisSchiela:2023:1} is compatible with our method.

The implementation of a practical algorithm that combines these two aspects is subject to future research.
Also the combination of our scheme with momentum based methods and/or variance reduction methods is a line of research that we plan to pursue in the future.
\nocite{Krishnamoorthy:2006:1}

\section*{Acknowledgments}

The authors used AI tools to enhance the written text regarding grammar and spelling, primarily in \cref{section:introduction}.

%% file: appendix.tex
\section{Proofs for Variance Bounds  and Asymptotic Behavior}
\label{section:proof-for-variance-behavior}

\lemmafirstboundonvariance*
\begin{proof}
	First, we compute:
	\begin{align*}
		\MoveEqLeft
		\norm{ f'_{\xi}(w) -  F'(w)}_{X^*}^2
		\\
		&
		=
		\norm{ f'_{\xi}(w) -  f'_{\xi}(w^\star)}_{X^*}^2
		+
		2 \, \inner{ f'_{\xi}(w) -  f'_{\xi}(w^\star)}{ f'_{\xi}(w^\star) -  F'(w^\star)}_{X^*}
		\\
		&
		\quad
		+ 2 \, \inner{ f'_{\xi}(w) -  f'_{\xi}(w^\star)}{F'(w^\star) -  F'(w)}_{X^*} + \norm{ f'_{\xi}(w^\star) -  F'(w^\star)}_{X^*}^2
		\\
		&
		\quad
		+ 2 \, \inner{ f'_{\xi}(w^\star) -  F'(w^\star)}{F'(w^\star) -  F'(w)}_{X^*} + \norm{ F'(w^\star) -  F'(w)}_{X^*}^2
		\\
		&
		\leq
		2 \, \norm{ f'_{\xi}(w) -  f'_{\xi}(w^\star)}_{X^*}^2
		+
		2 \, \norm{ f'_{\xi}(w^\star)}_{X^*}^2
		-
		2 \, \inner{f_\xi'(w) - f_\xi '(w^\star)}{F'(w)}_{X^*}
		\\
		&
		\quad
		- 2 \, \inner{ f'_{\xi}(w^\star)}{F'(w)}_{X^*}
		+
		\norm{ F'(w)}_{X^*}^2
		.
	\end{align*}
	Thus, taking the expectation yields:
	\begin{align*}
		\Var{\xi}{f'_\xi(w)}
		&
		\le
		2 \, \E[auto]{\xi}{L_\xi^2}\norm{w - w^\star}_X^2
		+
		2 \, \E[big]{\xi}{\norm{f'_\xi(w^\star)}_{X^*}^2}
		-
		\norm{F'(w)}_{X^*}^2
		\\
		&
		\le
		\paren[bigg](){2 \, \frac{\E{\xi}{L_\xi^2}}{\mu^2} - 1} \, \norm{F'(w)}_{X^*}^2
		+
		2 \, \E[big]{\xi}{\norm{f'_\xi(w^\star)}_{X^*}^2}
		.
	\end{align*}
\end{proof}

\lemmastrongerboundonvarianceuniformlyboundedLcase*
\begin{proof}
	As in the proof of \cref{lemma:first-bound-on-variance}, we get
	\begin{equation*}
		\Var[auto]{\xi}{f_\xi(w)}
		\le
		2 \, \E[auto]{\xi}{\norm{f_\xi'(w) - f_\xi'(w^\star)}_{X^*}^2}
		+
		2 \, \E[auto]{\xi}{\norm{f_\xi'(w^\star)}_{X^*}^2}
		-
		\norm{F'(w)}_{X^*}^2
		.
	\end{equation*}
	In the case that $L_\xi \le \Lmax$ holds uniformly in $\xi$, the authors in \cite{JohnsonZhang:2013:1} establish the bound\footnote{As in the proof of \cref{lemma:adapted-variance-bound} the setting of \cite{JohnsonZhang:2013:1} can easily be adapted to ours.}
	\begin{equation*}
		\E[big]{\xi}{\norm{f'_\xi(w) - f'_\xi(w^\star)}_X^2}
		\le
		2 \, \Lmax \, (F(w) - F(w^\star))
		.
	\end{equation*}
	Strong convexity yields $2 \, \mu \, (F(w) - F(w^\star)) \le \norm{F'(w)}_{X^*}^2$, and thus we obtain:
	\begin{equation*}
		\E[big]{\xi}{\norm{f'_\xi(w) - f'_\xi(w^\star)}_X^2}
		\le
		\frac{\Lmax}{\mu} \norm{F'(w)}_{X^*}^2
		.
	\end{equation*}
	We conclude for the variance:
	\begin{align*}
		\Var{\xi}{f'_\xi(w)}
		&
		\le
		2 \, \paren[Big](){\frac{\Lmax}{\mu} - 1} \norm{F'(w)}_{X^*}^2
		+
		2 \, \E[auto]{\xi}{\norm{f_\xi'(w^\star)}_{X^*}^2}
	\end{align*}
\end{proof}

\propositionunboundedvariationboundP*
\begin{proof}
	For $\gamma > 0$ and $\beta > 2$, let $\text{Par}(\gamma, \beta)$ be the Pareto distribution with parameters~$\gamma$ and~$\beta$.
	When $\xi \sim \text{Par}(\gamma, \beta)$, then we have $\E{\xi}{\xi} = \gamma \frac{\beta}{\beta - 1}$ and $\Var{\xi}{\xi} = \gamma^2 \frac{\beta}{(\beta - 2)(\beta - 1)^2}$; see, \eg, \cite[Chapter~23]{Krishnamoorthy:2006:1}.
	Thus with the choice $A_\xi \coloneqq \begin{psmallmatrix} \xi & 0 \\ 0 & 1 \end{psmallmatrix}$ and $f_\xi(w) \coloneqq \frac{1}{2} w^\transp A_\xi w$, $f_\xi$ is $\mu_\xi$-strongly convex and $L_\xi$-smooth with $\mu_\xi = \min(\xi, \mu)$, and $L_\xi = \max(\xi, 1)$.
	By definition, we have
	\begin{equation*}
		F(w)
		=
		\frac{1}{2} w^\transp
		\begin{pmatrix}
			\gamma \frac{\beta}{\beta - 1} & 0 \\ 0 & 1
		\end{pmatrix}
		w
		.
	\end{equation*}
	Suppose that $\eps > 0$ is arbitrary.
	When selecting $\gamma = \mu $ and $\beta = 2 + \eps$, then $F$ becomes $(\mu, 1)$-feasible for $\eps$ sufficiently small, and thus the corresponding stochastic optimization problem is $(\mu, 1)$-feasible.

	Further, choosing $w = \begin{psmallmatrix} s \\ 0	\end{psmallmatrix}$ with some scaling parameter~$s$, we observe
	\begin{itemize}
		\item
			$\Var{\xi}{\nabla f_\xi(w)} = s^2 \mu^2 \frac{\beta}{(\beta - 2)(\beta -1)^2} = s^2 \mu^2 \frac{2 + \eps}{\eps(1 + \eps)}$,

		\item
			$\norm{\nabla F(w)}_X^2 = s^2 \mu^2 \frac{\beta^2}{(\beta - 1)^2} = s^2 \mu^2 \paren[big](){\frac{2 + \eps}{1 + \eps}}^2$.
	\end{itemize}
	Thus, selecting $s = 2 \, \frac{\eps(1 + \eps)}{\mu V_0 (2 + \eps)}$, we obtain $w$ with $\Var{\xi}{\nabla F_\xi(w)} > V_0$ and
	\begin{equation*}
		\frac{\Var{\xi}{\nabla f_\xi(w)} - V_0}{\norm{\nabla F(w)}_X^2}
		\ge
		\frac{1}{4 \, \eps}
		.
	\end{equation*}
	Since $\eps$ was arbitrary, this proves the result.
\end{proof}

\propositionunboundedvariationboundPstar*
\begin{proof}
	This result is proved by a family of stochastic optimization problems that are strongly $(\mu, 1)$-feasible and satisfy
	\begin{equation*}
		\inf_{V_0 \in \R} V_1(V_0)
		\ge
		\frac{\paren[auto](){1-\mu}^3}{2 \, \mu \paren[auto](){2-\mu}^2}
		.
	\end{equation*}
	Then, using $\mu \le \frac{1}{2}$ gives the result.
	For $\mu \in (0,\frac{1}{2}]$ and $\alpha \in (0, 2 \pi)$ let
	\begin{equation*}
		A_1
		\coloneqq
		\begin{pmatrix}
			\mu \cos^2{\paren[auto](){\alpha }} + \sin^2{\paren[auto](){\alpha }}
			&
			\frac{\paren[auto](){1 - \mu} \sin{\paren[auto](){2 \alpha }}}{2}
			\\
			\frac{\paren[auto](){1 - \mu} \sin{\paren[auto](){2 \alpha }}}{2}
			& \
			\mu \sin^2{\paren[auto](){\alpha }} + \cos^2{\paren[auto](){\alpha }}
		\end{pmatrix}
	\end{equation*}
	and
	\begin{equation*}
		A_2
		\coloneqq
		\begin{pmatrix}
			\mu \cos^2{\paren[auto](){\alpha }} + \sin^2{\paren[auto](){\alpha }}
			&
			-\frac{\paren[auto](){1-\mu} \sin{\paren[auto](){2 \alpha }}}{2}
			\\
			-\frac{\paren[auto](){1-\mu} \sin{\paren[auto](){2 \alpha }}}{2}
			&
			\mu \sin^2{\paren[auto](){\alpha }} + \cos^2{\paren[auto](){\alpha }}
		\end{pmatrix}
		.
	\end{equation*}
	Then,
	\begin{equation*}
		A
		=
		\frac{1}{2}(A_1 + A_2)
		=
		\begin{pmatrix}
			\mu \cos^2{\paren[auto](){\alpha }} + \sin^2{\paren[auto](){\alpha }}
			&
			0
			\\
			0
			&
			\mu \sin^2{\paren[auto](){\alpha }} + \cos^2{\paren[auto](){\alpha }}
		\end{pmatrix}
		.
	\end{equation*}
	As is easily checked, $A_1$ and $A_2$ have the eigenvalues~$\mu$ and~$1$, and~$A$ has the eigenvalues $\mu \cos^2{\paren[auto](){\alpha}} + \sin^2{\paren[auto](){\alpha }}$ and $\mu \sin^2{\paren[auto](){\alpha }} + \cos^2{\paren[auto](){\alpha }}$.

	For $i = 1, 2$, let $f_i(w) \coloneqq \frac{1}{2} w^\transp A_1 w$ and $F(w) \coloneqq \frac{1}{2}(f_i(w) + f_2(w))$.
	The corresponding SOP (with $X = \{1,2\}$ and $P$~being the uniform distribution on~$X$) is strongly $(\mu, L)$-feasible.
	Trivially, $\nabla f_i(w) = A_i w$ and $\nabla F(w) = A w$.
	When fixing $\alpha = \arcsin \paren[normal](){\sqrt{\mu}}$ and choosing $w = s \begin{psmallmatrix} \frac{1}{\mu} \\ 0	\end{psmallmatrix}$ for a scaling parameter $s > 0$, we observe
	\begin{itemize}
		\item
			$\Var{\xi}{\nabla f_\xi(w)} = {\frac{s^2 \paren[auto](){1-\mu}^3}{\mu}}$
			,

		\item
			$\norm{\nabla F(w)}_X^2 = s^2 \paren[auto](){\mu - 2}^2$
			.
	\end{itemize}
	Thus, choosing $s \coloneqq \sqrt{2} \sqrt{\frac{V_0 \mu}{\paren[auto](){1 - \mu}^3}}$ provides us with a vector~$w$ with $\Var{\xi}{\nabla f_\xi(w)} > V_0$.
	Therefore,
	\begin{equation*}
		V_1(V_0)
		\ge
		\frac{\Var{\xi}{\nabla f_\xi(w)} - V_0}{\norm{\nabla F(w)}_X^2}
		=
		\frac{\paren[auto](){1-\mu}^3}{2 \, \mu \paren[auto](){2-\mu}^2}
		.
	\end{equation*}
\end{proof}

\begin{remark}
	\label{remark:geometric-motivation-for-variance-result}
	The example given in the proof of \cref{proposition:unbounded-variation-bound-Pstar} has a clear geometric motivation.
	The matrices $A_i$ have the form $Q_i^\transp D Q_i$, where
	\begin{equation*}
		Q_i
		=
		\begin{pmatrix}
			(-1)^{i +1} \cos{\paren[auto](){\alpha}}
			&
			- \sin{\paren[auto](){\alpha }}
			\\
			\sin{\paren[auto](){\alpha }}
			&
			(-1)^{i +1} \cos{\paren[auto](){\alpha}}
		\end{pmatrix}
	\end{equation*}
	is a rotation matrix with angle $(-1)^{i+1} \alpha$, and $D = \begin{psmallmatrix} \mu & 0 \\ 0 & 1 \end{psmallmatrix}$.
	Then, $A_i$ clearly has the eigenvalues $\mu$ and $1$, and $\begin{psmallmatrix} (-1)^{i+1} \cos(\alpha) \\ \sin(\alpha)	\end{psmallmatrix}$ is an eigenvector for the eigenvalue~$\mu$ and $\begin{psmallmatrix} -\sin(\alpha) \\ (-1)^{i + 1} \cos(\alpha) \end{psmallmatrix}$ is an eigenvector for the eigenvalue~$1$.
	An eigenvector for the smallest eigenvalue $\mu \cos^2(\alpha) + \sin^2(\alpha)$ of~$A$ is given by~$w$ as selected in the proof.
	By choosing~$w$ in this way, $\norm{\nabla F(w)}_X$ is only affected by the smallest eigenvalue, here converging to zero, while the variance $\Var{\xi}{\nabla f_\xi(w)}$ is also affected by the larger eigenvalue~$1$.
	By sending $w$~to infinity in the right way (as described above, taking into account that $V_0$ can be arbitrary), one takes advantage of this fact and uses that $\nabla F(w)$ scales linearly with~$\mu$, while this is not the case for the variance.

	The result further shows that the bound in \cref{lemma:stronger-bound-on-variance-uniformly-bounded-L-case} is asymptotically sharp.
\end{remark}

\section{Complete Algorithm in Pseudo Code}
\label{section:complete-algorithm-in-pseudo-code}

Based on the results in this paper, we propose an extension of SGD, which incorporates our estimation techniques.
To write our algorithm in clear pseudo code, we first introduce the following well known algorithms for line search and exponential smoothing.
\Cref{algo:line search} defines a simple line search, used for initialization and in the case of unsuccessful steps.
In the line search procedure, we select a new realization of $\xi$ in each iteration and thereby reduce the influence of one single batch to the overall training process.
This prevents unnecessarily small step sizes due to outliers.

Based on the concept presented in \cref{subseciton:exponential-smoothing}, \cref{algo:smoothing} defines the averaging technique imposed to average over individual observations and thus to obtain a stable estimator.
It is based on the estimation and smoothing techniques explained in \cref{section:estimation-techniques}.

Finally, in \cref{algorithm:adaptive-SGD}, we propose the adaptive SGD scheme in more detail.

\begin{algorithm}[htp]
	\caption{$\SLS((f_\xi, \Omega, P), H, \alpha, w, \eta)$}
	\label{algo:line search}
	\begin{algorithmic}[1]
		\Require
		SOP $(f_\xi, \Omega, P)$,
		preconditioner $H$,
		initial step size $\alpha > 0$,
		current position~$w$,
		shrinking parameter $\eta < 1$

		\State Sample $\xi \sim P$
		\While{$f_\xi(w - \alpha \, \inv{H} f_\xi'(w)) > f_\xi(w)$}
		\State $\alpha \gets \eta \, \alpha$
		\State Sample new $\xi \sim P$
		\EndWhile
		\State \Return $\alpha$, $w - \alpha \, \inv{H} f_\xi'(w)$, $\xi$
	\end{algorithmic}
\end{algorithm}

\begin{algorithm}[htp]
	\caption{$\ExpSmooth(q, q^+, \gamma)$}
	\label{algo:smoothing}
	\begin{algorithmic}[1]
		\Require
		current value $q$,
		next observation $q^+$,
		discount factor $\gamma \in (0,1)$
		\State $q \gets \gamma \, q + (1 - \gamma) \, q^+$
		\State \Return $q$
	\end{algorithmic}
\end{algorithm}

\begin{algorithm}[htp]
	\caption{SGD with adaptive step size control}
	\label{algorithm:adaptive-SGD}
	\begin{algorithmic}[1]
		\Require
		pointwise $(\mu_\xi, L_\xi)$-feasible SOP $(f_\xi, \Omega, P)$,
		line search parameter $\eta_\alpha < 1$,
		initial step size $\alpha$,
		preconditioner $H$,
		discount factors $\gamma_k$
		\State $\alpha_0, w_1, \xi_0 \gets \SLS((f_\xi, \Omega, P), H, \alpha, w_0, \eta_\alpha)$.
		\State $L \gets \frac{1}{\alpha_0}$
		\State $\variance \gets 0$
		\State $\gradnormsquaredaveraged \gets \norm{f'_{\xi_0}(w_0)}_{X^*}^2$
		\State $\evallinesearch \gets \False$
		\For{$k \ge 1$}
		\If{\evallinesearch}
		\State $\alpha_k, w_{k+1}, \xi_k \gets \SLS((f_\xi, \Omega, P), H, \alpha_k, w_k, \eta_\alpha)$
		\State $\evallinesearch \gets \False$
		\Else
		\State Sample $\xi_k \sim P$
		\State $w_{k+1} \gets w_k - \alpha_k \inv{H} f'_{\xi_k}(w_k)$
		\If{$f_{\xi_k}(w_{k+1}) > f_{\xi_k}(w_k)$}
		\State $w_{k+1} \gets w_k$
		\Comment{reject unsuccessful step}
		\State $\evallinesearch \gets \True$
		\Comment{perform line search in next iteration}
		\State Continue with next~$k$
		\EndIf
		\EndIf
		\State $\displaystyle \overline{\variance} \gets \frac{f_{\xi_k}(w_k) - f_{\xi_{k-1}}(w_k)}{\alpha_{k-1}}$
		\State $\variance \gets \ExpSmooth(\variance, \overline{\variance}, \gamma_k)$
		\State $\displaystyle \overline L \gets 2 \, \frac{f_{\xi_k}(w_{k+1}) - f_{\xi_k}(w_k) + \alpha \, \norm{f_{\xi_k}(w_k)}_{X^*}^2}{\alpha^2 \norm{f_{\xi_k}'(w_k)}_{X^*}^2}$
		\State $L \gets \ExpSmooth(L, \overline L, \gamma_k)$
		\State $\gradnormsquaredaveraged \gets \ExpSmooth(\gradnormsquaredaveraged, \norm{f'_{\xi_k}(w_k)}_{X^*}^2, \gamma_k)$
		\State $\displaystyle \overline \alpha \gets \frac{\gradnormsquaredaveraged - \variance}{L \, \gradnormsquaredaveraged}$
		\State $\alpha \gets \ExpSmooth(\alpha, \overline \alpha, \gamma_k)$
		\EndFor
	\end{algorithmic}
\end{algorithm}
\FloatBarrier